\documentclass[letterpaper, paper,10pt]{AAS}		

\usepackage{bm}
\usepackage{amsmath}
\usepackage[colorlinks=true, pdfstartview=FitV, linkcolor=black, citecolor= black, urlcolor= black]{hyperref}
\usepackage{overcite}
\usepackage{footnpag}			      	
\usepackage[ruled,linesnumbered]{algorithm2e}
\SetKwIF{If}{ElseIf}{Else}{if}{then}{else if}{else}{endif}
\usepackage{graphicx}
\usepackage{subfig}

\PaperNumber{22-581}

\begin{document}

\title{CLOSED-FORM OPTIMAL SOLUTIONS FOR PROPULSIVE-DIFFERENTIAL DRAG CONTROL OF SPACECRAFT SWARMS}

\author{Matthew Hunter\thanks{Ph.D. Candidate, Stanford University, Department of Aeronautics \& Astronautics, Space Rendezvous Laboratory, Durand Building, 496 Lomita Mall Stanford, CA 94305}
\ and Simone D'Amico\thanks{Assistant Professor, Stanford University, Department of Aeronautics \& Astronautics, Space Rendezvous Laboratory, Durand Building, 496 Lomita Mall Stanford, CA 94305}
}

\maketitle{}

\begin{abstract}
This paper presents a novel optimal control approach for Distributed Space Systems (DSS) to minimize the total delta-v cost for large, time-constrained reconfigurations. Reconfiguration cost is decreased by blending differential drag and propulsive control within the same control window. By parameterizing the problem with Relative Orbital Elements (ROE) and leveraging reachable set theory, a maneuver planning methodology is derived to produce provably-optimal, full roto-translational control solutions in closed-form, such that they can be implemented on hardware-limited spacecraft.
\end{abstract}

\section{Introduction}
Distributed Space Systems (DSS) present new opportunities for space missions by using multiple spacecraft to accomplish objectives that would be hard to realize by a single monolithic spacecraft. A desire for autonomous operation and the trend towards smaller spacecraft place two major constraints on DSS mission capability: individual spacecraft fuel budget and on-board processing power. These motivate interest in closed-form fuel-optimal guidance and control algorithms that are computationally tractable, easy to test, and robust in a variety of orbit regimes and perturbative conditions. Rather than treating perturbative forces as estimation and control errors, additional delta-v savings can be found by incorporating environmental forces dependent on differential spacecraft attitude into the control problem. This is especially true in Low Earth Orbit (LEO), where atmospheric drag is the largest perturbation force behind J2, for missions that require large changes in along-track separation to accomplish science objectives, such as the planned Space Weather Atmospheric Reconfigurable Multiscale Experiment (SWARM-EX) mission through the National Science Foundation (NSF). \cite{swarm-ex}

Previous research has investigated the use of differential drag in LEO to conduct reconfigurations or formation keeping with little or no delta-v cost. Differential drag suffers from a lack of control authority due to uncertainty and variability associated with several of its critical parameters, extremely limited instantaneous control magnitude, and restricted directionality. Therefore, the literature primarily focuses on robust control through iterative techniques to converge to a desired final state using no thrust. This is especially common for maintaining passive safety for DSS with many individuals. Riano-Rios et al. integrated a Linear Quadratic Regulator (LQR) with constrained least squares to provide simultaneous drag maneuvers for multiple deputy spacecraft to converge to a desired overall formation and determined that adding an intermediary along-orbit formation greatly reduced the chance of collision. \cite{bev_drag_leastsquares} Mazal et al. was able to mathematically guarantee convergence with an LQR controller by placing bounds on the final state accuracy, given known uncertainty in the atmospheric density and spacecraft drag coefficient. \cite{bev_drag_rendezvous} Traub et al. emphasized the preservation of mechanical energy and mission lifetime by solving a maneuver planning optimization problem that minimizes the loss of semi-major axis in a spacecraft's orbit through Model Predictive Control (MPC). \cite{traub_drag} Drag-only formation maintenance has also been proven to be effective for large satellite constellations, such as the Planet Flock. \cite{planet} 


Parallel efforts have created more comprehensive roto-translational maneuver schemes that simultaneously stabilize the rotational dynamics, or attitude, of the spacecraft while still accomplishing translational control, or changing the orbit of the spacecraft. Advances in dual-quaternion formulations have allowed iterative control of six degrees of freedom (6-DOF) dynamic systems with provable Lyapunov stability, but this approach struggles to characterize perturbations and has yet to produce a provably optimal control law. \cite{stanfield_dualq_lqr,stanfield_dualq_apps,yang_dualq} Riano-Rios et al. accomplishes roto-translational Lyapunov-stability for the narrow subset of spacecraft equipped with Drag Maneuvering Devices, which allow for partially independent control of translational forces and rotational torques due to drag. \cite{bev_drag_roto-trans} Pastorelli et al. used a similar approach to achieve 6-DOF Lyapunov stability by modulating a spacecraft's center of pressure with a drag sail placed on the back-end of the spacecraft. \cite{bev_drag_propellantless}

Very little work in differential drag control has taken advantaged of the benefits offered by Integration Constant (IC) states, such as Relative Orbital Elements (ROE), to reduce nonlinear relative dynamics to linear models. In the broader sense of extended or low thrust maneuvers, Di Mauro et al. examined non-impulsive maneuvers in near-circular orbit to produce maneuver schemes in closed-form for piecewise continuous maneuvers and numerically by posing fuel-optimal control as a mixed-integer linear program. \cite{bev_extendedman,bev_finite-time} Scala et al. relies on ROE's to provide a linear dynamic model to produce affine constraints for a disciplined convex programming problem for low-thrust control. \cite{scala_lowthrust} Ben-Larbi et al. also took advantage of linearized ROE dynamics to produce a feedforward and feedback Proportional Integral Derivative (PID) control law to conduct input-output linearization of the complex nonlinear variations of drag control magnitude. \cite{ben-larbi_drag_ROE} Another ROE approach by Koenig et al. models differential drag as on-board low-thrust control and uses switching lanes to control about a desired state. \cite{drag_koenig} However, none of these techniques consider both propulsive and drag control in the same control window, for situations when drag control alone will not reach the desired state in a given time. Furthermore, drag and propulsive control complement one another in simultaneous application, as drag control contributes maneuvering without fuel consumption while propulsive control provides the controllability that drag lacks.

This work proposes a methodology to produce provably optimal closed-form solutions for hybrid differential drag - propulsive control that provide delta-v cost savings for large, time-constrained reconfigurations over propulsive-only approaches, while maintaining full controllability. Realizing hybrid control in closed form requires two major contributions to the state of the art. First, the hybrid control problem is formulated as a linear time-variant (LTV) system in ROE-space by treating differential drag as a control force and analytically integrating its effect. Second, this convex formulation is characterized through reachable set theory, decoupling the hybrid control problem into successive drag and propulsive maneuver planning steps without compromising optimality to allow optimal drag profiles to be found in closed form. This paper presents the methodology used to reach these solutions in four sections. First, the LTV hybrid control problem is derived. Second, reachable set theory is used to visualize the controllability of differential drag and separate optimal drag and propulsive control solutions. Third, all possible optimal drag profile variations are found analytically. Fourth, a flexible drag profile architecture is designed to bridge the gap between possible optimal drag profiles, as well as the gap between hybrid and drag-only control. Finally, the validity of the proposed closed-form methodology is proven by comparing its solutions with those produced from numerical interior point methods, and the validity of the hybrid approach overall is verified in a high fidelity orbital propagator that integrates the Gauss Variational Equations (GVE) subject to relevant perturbations.

\section{Background}

\subsection{Relative Orbital Elements as State Parameters}
Spacecraft state estimation and control actions are dictated in cartesian frames such as Earth-Centered-Inertial (ECI) or Radial-Tangential-Normal (RTN). ECI is a geocentric right-handed triad aligned with the equator, vernal equinox, and North Pole. RTN centers on a spacecraft's center of mass and aligns with the orbital plane, positive radially outward and tangentially along the direction of travel, and the angular momentum vector. However, more recent control efforts have used intermediary Integration Constant (IC) states to build accurate linear models of nonlinear systems for analytical control algorithms and guidance strategies. A commonly used IC parameterization is Relative Orbital Elements (ROE), which were pioneered and flight proven on the PRISMA and TanDEM-X missions to demonstrate deterministic propulsive maneuver schemes for reconfiguration and formation keeping. \cite{damico_phd} ROE are defined as a nonlinear combination of the Keplerian orbital elements and therefore are slowly varying under the presence of both Keplerian and perturbation dynamics. Further osculating dynamics can also be eliminated by adopting mean ROE from mean orbital elements, which approximate a spacecraft's trajectory as an ideal ellipse by removing perturbative effects from the osculating orbital elements. For Keplerian motion and small separations compared to orbit radius, ROE are equivalent to the IC's of the linear differential equations of relative motion, such as the Hill-Clohessy-Wiltshire (HCW) and Yamanaka-Ankersen (YA) equations. This work adopts the following definition of ROE
\begin{equation} \label{eq: ROE}
    \delta\boldsymbol{\alpha} = 
    \begin{bmatrix}
    \delta a \\ \delta \lambda \\ \delta e_x \\ \delta e_y \\ \delta i_x \\ \delta i_y
    \end{bmatrix} = 
    \begin{bmatrix}
    \frac{a_d-a_c}{a_c} \\
    (M_d-M_c) + (\omega_d-\omega_c) + \cos(i_c)(\Omega_d-\Omega_c) \\
    e_d\cos(\omega_d) - e_c\cos(\omega_c) \\
    e_d\sin(\omega_d) - e_c\sin(\omega_c) \\
    i_d - i_c \\
    \sin(i_c)(\Omega_d-\Omega_c)
    \end{bmatrix}
\end{equation}
\noindent where $\delta a$ is the mean relative semi-major axis, $\delta \lambda$ is the mean relative longitude, $\delta e_x$ and $\delta e_y$ are the components of the mean relative eccentricity vector $\delta \boldsymbol{e}$, $\delta i_x$ and $\delta i_y$ are the components of the mean relative inclination vector $\delta \boldsymbol{i}$, and $a$, $e$, $i$, $\Omega$, $\omega$, and $M$ are the mean Keplerian orbital elements. These ROE are referred to as quasi-nonsingular because they are valid for inclined orbits of arbitrary eccentricity. This formulation describes the relative dynamics between a chief and deputy spacecraft, denoted by subscripts $c$ and $d$ respectively. The chief can describe a reference satellite that does not conduct maneuvers or a reference orbit. Isolating the long term and secular dynamics enables the derivation of accurate ROE-based State Transition Matrices (STM) to represent relative spacecraft motion as a LTV system. Further advances in STM’s by Koenig and Guffanti have included corrections for the Earth's oblateness (J2), atmospheric drag, solar radiation pressure (SRP), and third body Sun and Moon effects in closed-form, further increasing the accuracy of the model over long time steps. \cite{STM_koenig,STM_guffantiSRP,STM_guffanti3rd} In the context of maneuver planning, STM's propagate the effect of control forces from the time of actuation to the end of the control window. For notational simplicity, auxiliary matrix $\Gamma$ is introduced to quantify the change in mean ROE final state at the end of the control window $t_f$ given an instantaneous RTN control action at time $t$ and is given as
\begin{align} \label{eq: Gamma matrix}
    &\boldsymbol{\Gamma}(t)=\boldsymbol{\Gamma}(\boldsymbol{\alpha}_c(t),t_f-t) = \boldsymbol{\Phi}^{J_2} (\boldsymbol{\alpha}_c(t),t_f-t) \boldsymbol{B}(\boldsymbol{\alpha}_c(t)) \\
    & = 
    \begin{bmatrix}
    1 & 0 & 0 & 0 & 0 & 0 \\
    \Phi_{21}^{J_2} & 1 & 0 & 0 & \Phi_{25}^{J_2} & 0 \\
    0 & 0 & \cos(\Dot{\omega}_c(t_f-t)) & -\sin(\Dot{\omega}_c(t_f-t)) & 0 & 0 \\
    0 & 0 & \sin(\Dot{\omega}_c(t_f-t)) & \cos(\Dot{\omega}_c(t_f-t)) & 0 & 0 \\
    0 & 0 & 0 & 0 & 1 & 0 \\
    \Phi_{61}^{J_2} & 0 & 0 & 0 & \Phi_{65}^{J_2} & 1 \\
    \end{bmatrix} 
    \frac{1}{a_c n_c}
    \begin{bmatrix}
    0 & 2 & 0 \\
    -2 & 0 & 0 \\
    \sin(u_c) & 2\cos(u_c) & 0 \\
    -\cos(u_c) & 2\sin(u_c) & 0 \\
    0 & 0 & \cos(u_c) \\
    0 & 0 & \sin(u_c) \\
    \end{bmatrix} \nonumber
\end{align}
\noindent with the following simplifying substitutions
\begin{align}
    \eta_c = \sqrt{1-e_{c}^{2}} \text{ , } \kappa_c = \frac{3}{4}\frac{J_2 R^{2}_{E} \sqrt{\mu}}{a_{c}^{7/2} \eta_{c}^{4}} \text{ , } \Dot{\omega}_c = \kappa_c(5\cos^2&(i_c)-1) \\
    \Phi_{21}^{J_2} = - \left( \frac{3}{2}n_c+\frac{7}{2}\kappa_c(1+\eta_c)(3\cos^2(i_c)-1) \right)& (t_f-t) \nonumber \\
    \Phi_{25}^{J_2} = -\kappa_c(4+3\eta_c)\sin(2i_c)(t_f-t) \text{ , } \Phi_{61}^{J_2} = \frac{7}{2}\kappa_c\sin(2i_c)(t_f-t) \text{ , } &\Phi_{65}^{J_2} = 2\kappa_c\sin^2(i_c)(t_f-t) \nonumber
\end{align}
\noindent where STM $\boldsymbol{\Phi}^{J_2}$ includes the Koenig corrections for J2 for mean ROE in near-circular orbit, the control input matrix $\boldsymbol{B}$ maps control actions in the RTN frame to an equivalent change in osculating or mean ROE in near-circular orbit, $\boldsymbol{\alpha}_c$ is the chief mean orbital elements at the time of the control action, $n_c$ is the mean motion, $\Dot{\omega}_c$ denotes the rate of drift of the argument of periapsis, and $u_c$ is the mean argument of latitude at the time of the control action. This dynamics model assumes that the spacecraft swarm is operating in LEO where the dominant perturbations to Keplerian orbit are J2 and atmospheric drag, and orbits can be approximated as near-circular. This is a reasonable assumption, as the magnitude of differential drag is higher in LEO than any other orbit regime, making this regime the most relevant for this work. Given that atmospheric drag will be posed as a form of control in later sections, the only relevant correction required for state propagation accuracy over long control windows is the J2 correction. Additionally, an important consequence of the ROE linearization process is the decoupling of the in-plane ROE dynamics ($\delta a$, $\delta \lambda$, $\delta \boldsymbol{e}$) from the out-of-plane ROE dynamics ($\delta \boldsymbol{i}$) for this choice of ROE in near-circular orbit.

\subsection{Reachable Set Theory for Impulsive Control}
Reachable set theory provides a geometric framework for optimal control problems with norm-like cost functions and LTV dynamics. The framework consists of contact and support functions that can be derived as analytical expression of the spacecraft's state. Previous efforts by Koenig, Chernick, and D'Amico used reachable set theory to derive optimality conditions to produce efficient, globally optimal solvers for impulsive control of LTV dynamics, and the methodology of this approach is summarized here. \cite{chernick,impulse_koenig} Let $U(c,t)$ be the set of control actions $\boldsymbol{u}$ in the RTN frame with a cost no greater than $c$ (L2 norm cost function) at time $t$. Let $S(c,t)$ be the set of pseudostates $\Delta \delta \boldsymbol{\alpha}$ that can be reached by a single control input such that $\boldsymbol{u} \in U(c,t)$, given as
\begin{equation} \label{eq: RST single maneuver}
    S(c,t) = \{\Delta \delta \boldsymbol{\alpha}: \Delta \delta \boldsymbol{\alpha} = \boldsymbol{\Gamma}(\boldsymbol{\alpha}_c(t),t_f-t)\boldsymbol{u}, \boldsymbol{u} \in U(c,t)\}
\end{equation}
\noindent where $\Delta \delta \boldsymbol{\alpha} = \delta\boldsymbol{\alpha}(t_f) - \boldsymbol{\Phi}(\boldsymbol{\alpha}_{c}(t_{0}),t_{f}-t_{0})\delta\boldsymbol{\alpha}(t_{0})$ is invariant with respect to the control problem, $\boldsymbol{\Gamma}$ is given in Eq. \ref{eq: Gamma matrix}, and $t_0$ and $t_f$ define the beginning and end of the control window. $\Delta \delta \boldsymbol{\alpha}$ and $\boldsymbol{\Gamma}$ are substituted to simplify notation. Let $T$ be a discrete set of maneuver times in the interval $[t_0,t_f]$, and let $S(c,T)$ be the set of $\Delta \delta \boldsymbol{\alpha}$ that can be reached using a single control action applied at any time in $T$ with a cost no greater than $c$, given as
\begin{equation}
    S(c,T) = \bigcup_{t\in T} S(c,t)
\end{equation}
\noindent Finally, let $S^*(c,T)$ be the set of $\Delta \delta \boldsymbol{\alpha}$ that can be reached by a set of $k \geq 1$ control inputs applied at times in $T$ with a combined cost no greater than $c$, defined as
\begin{equation} \label{eq: RST convex hull}
    S^*(c,T) = \{\Delta \delta \boldsymbol{\alpha}: \Delta \delta \boldsymbol{\alpha} = \sum^{k}_{j=1} \boldsymbol{\Gamma}(\boldsymbol{\alpha}_c(t_j),t_f-t_j)\boldsymbol{u}_j, t_j \in T, \boldsymbol{u}_j \in U(c_j,t_j), \sum^{k}_{j=1} c_j = c \}
\end{equation}
$\boldsymbol{\Gamma}(\boldsymbol{\alpha}_c(t_j),t_f-t_j)\boldsymbol{u}_j$ must be in $S(c_j,T)$ for any $t_j \in T$ and $\boldsymbol{u}_j \in U(c_j,t_j)$. Because the cost of a control input scales linearly with its magnitude,  $S^*(c,T)$ forms the convex hull of $S(c,T)$, and the minimum cost for a desired reconfiguration must lie on the boundary of $S^*(c,T)$. Chernick proved that a 2$n$-dimensional (2$n$D) control problem can be decomposed into $n$ 2D planes and evaluated separately to find the minimum cost of the reconfiguration without loss of generality. \cite{chernick} The coupling between planes instead results in optimality constraints, such as optimal maneuver times, which establish a system of equations for an optimal closed-form solution.

For a 6-dimensional ROE state, a convenient choice of 2D planes is the $\Delta \delta \boldsymbol{a}$ plane ($\Delta \delta a$, $\Delta \delta \lambda$), the $\Delta \delta \boldsymbol{e}$ plane, and the $\Delta \delta \boldsymbol{i}$ plane, where the leading $\Delta$ denotes a pseudostate as defined in Eq. \ref{eq: RST single maneuver}. The dominant plane of a reconfiguration refers to the plane that requires the highest cost to reach the desired pseudostate in its plane, and the dominant ROE, or dominance case, indicates the specific ROE that requires the highest cost to reach its desired pseudostate within the dominant plane. When applying this methodology to spacecraft, the cost associated with a specific dominance case is known as its delta-v minimum, the highest of which drives the overall delta-v minimum of the reconfiguration, and Chernick derived an analytical formulation of the delta-v minimum for each ROE dominance case.\cite{chernick} Because ROE linearization decouples in-plane and out-of-plane dynamics, in-plane reconfigurations can be analyzed and conducted separately from out-of-plane reconfigurations, turning a 6D problem into a set of equivalent 4D and 2D sub-problems.

\section{Problem Definition}
The hybrid control problem proposed below expands on the propulsive control formulation by Chernick and D’Amico without violating any of the assumptions or key characteristics required to produce closed-form solutions. Using the LTV dynamics model described previously, the most general formulation of the reachable set theory-compatible optimal control problem is
\begin{align} \label{eq:general control problem}
    \text{minimize  } c =& \int_{T} f(\boldsymbol{u}(t),t) \\
    \text{subject to  } \delta\boldsymbol{\alpha}(t_{f}) - \boldsymbol{\Phi}(\boldsymbol{\alpha}_{c}(t_{0}),&t_{f}-t_{0})\delta\boldsymbol{\alpha}(t_{0}) = \Delta \delta \boldsymbol{\alpha} = 
    \int_{T} \boldsymbol{\Gamma}(\boldsymbol{\alpha}_c(t),t_f-t)  \boldsymbol{u}(t) dt \nonumber
\end{align}
\noindent where $\delta\boldsymbol{\alpha}(t_{0})$ describes the deputy spacecraft’s initial mean ROE state, $\delta\boldsymbol{\alpha}(t_{f})$ is the deputy spacecraft's final mean ROE state, $f$ is a norm-like cost function, $\boldsymbol{\Gamma}$ is given in Eq. \ref{eq: Gamma matrix}, and $T$ is a set of times that span the control window defined by the interval $[t_0, t_f]$. The difference between the final desired and the initial state propagated by the STM across the control window is again known as the pseudostate $\Delta \delta \boldsymbol{\alpha}$ and is invariant with respect to the control problem. The control itself can be split into two types, propulsive and drag, defined as
\begin{align} \label{eq:general expanded control problem}
    \text{minimize  } c =& \sum_{j=1}^{k} ||\boldsymbol{u}^{\text{propulsive}}_{j}||_{2} + \int_{T} f(\boldsymbol{u}^{\text{drag}}(t),t) \\
    \text{subject to  } \Delta \delta \boldsymbol{\alpha} =& 
    \sum_{j=1}^{k} \boldsymbol{\Gamma}(t_{j}) \boldsymbol{u}^{\text{propulsive}}_{j} + 
    \int_{T} \boldsymbol{\Gamma}(t)  \boldsymbol{u}^{\text{drag}}(t) dt \text{ , } t_j \in T \nonumber
\end{align}
Propulsive maneuvers are characterized by high thrust magnitude and short maneuver duration with respect to the orbit period, allowing them to be approximated as instantaneous additions of delta-v or “impulsive” maneuvers. Propulsive control can then be posed as an affine function of impulsive maneuvers and separated from the remaining integral of time-varying control. No assumptions about the nature of drag control are required to separate it from propulsive maneuvers, except that it can be conducted independently of propulsive control. 

\section{Differential Drag Formulation}
Unlike conservative forces such as the J2 perturbation, atmospheric drag can be controlled on-board each spacecraft through the modulation of attitude or cross-sectional area normal to the direction of travel. The force of atmospheric drag on a single spacecraft is given by
\begin{equation}
    p^{{\text{drag}}} = \frac{1}{2}\rho v^2 \frac{C_D A}{m} = \frac{1}{2}\rho n^2 a^2 B
\end{equation}
\noindent where $\rho$ is the atmospheric density, $v$ defines the spacecraft velocity, $C_D$ is the spacecraft's coefficient of drag, $A$ is the cross-sectional area, $m$ is the spacecraft's mass, $B = \frac{C_D A}{m}$ defines the spacecraft ballistic coefficient, and $n$ and $a$ are the mean motion and mean semi-major axis of the spacecraft's orbit. In near-circular orbit, $v = na$ approximates the spacecraft velocity as the tangential velocity for circular motion and is constant over the entire orbit. The force from atmospheric drag acts opposite the direction of travel and therefore always acts in the negative tangential direction. The force also has a minimum and maximum limit, dictated by the spacecraft's minimum and maximum achievable cross-sectional area and corresponding ballistic coefficient.

Differential drag is the difference between the force of atmospheric drag on the chief and deputy spacecraft and is controlled by modulating differential attitude between the spacecraft. This differential force is expressed in the RTN frame as
\begin{align} \label{eq:drag force}
    \boldsymbol{u}^{\text{drag}}(t) = 
    \begin{bmatrix}
    0 \\ p^{\text{drag, chief}}(t) \\ 0
    \end{bmatrix} &- 
    \begin{bmatrix}
    0 \\ p^{\text{drag, deputy}}(t) \\ 0
    \end{bmatrix} = 
    \begin{bmatrix}
    0 \\ \delta p^{\text{drag}}(t) \\ 0
    \end{bmatrix} = 
    \begin{bmatrix}
    0 \\ \frac{1}{2}n_c^2a_c^2\Delta B^{\rho}(t) \\ 0
    \end{bmatrix} 
\end{align}
\noindent where $\Delta B^{\rho}(t) = \rho_c(t)\Delta B_c(t) - \rho_d(t)\Delta B_d(t)$ defines the augmented differential ballistic coefficient for notational simplicity. The chief and deputy are assumed to be in close orbits, such that they both have the chief's mean motion and mean semi-major axis. Additionally, assuming relatively small separations between the chief and deputy allows the difference between the tangential directions of the spacecraft to be neglected, and the force of atmospheric drag acts only in the tangential direction with respect to the chief on both spacecraft. This results in a differential force that also only acts in the tangential direction and has minimum and maximum limits depending on the drag force limits achievable by the chief and deputy. This assumption is reasonable, as a spacecraft conducting relative maneuvers with respect to a chief to maintain a formation inherently implies an along-track separation that is significantly smaller than the orbital radius, or relative maneuvering would not be relevant.

To analyze differential drag with reachable set theory, drag needs to be posed in an affine LTV formulation. This can be achieved by discretizing the control window into a series of time steps and calculating the effect that a constant drag maneuver over a given time step will have on the pseudostate, formulated as 
\begin{equation} \label{eq: affine drag}
    \int_{T} \boldsymbol{\Gamma}(t)  \boldsymbol{u}^{\text{drag}}(t) dt = 
    \sum_{i=1}^{m} \int_{t_{i}}^{t_{i+1}} \boldsymbol{\Gamma}(t_{i})  \boldsymbol{u}^{\text{drag}}_{i}(t_{i}) dt = \sum_{i=1}^{m} \boldsymbol{\Gamma}^{\text{drag}}(t_{i}) \Delta B^{\rho}_{i}
\end{equation}
where the control window has been divided into $m$ time steps and the control variable can be posed as a vector of $\Delta B^{\rho}$ values, one value for each time step. Because differential drag is an extended maneuver, the effect of drag on the pseudostate requires integration with respect to time. The integration is of Eq. \ref{eq: affine drag} accomplished analytically and results in the following expression
\begin{align} \label{eq: integrated drag}
    &\int_{t_{i}}^{t_{i+1}} \boldsymbol{\Phi}^{J_2} (\boldsymbol{\alpha}_c(t),t_f-t) \boldsymbol{B}(\boldsymbol{\alpha}_c(t)) dt = \int_{t_{i}}^{t_{i+1}} \boldsymbol{\Gamma}(t)  \boldsymbol{u}^{\text{drag}}(t) dt \\
    &= a_c n_c
    \begin{bmatrix}
    t_{i+1}-t_i \\
    0.5\Phi_{\text{21}}^{\text{J2}}(t_i-t_{i+1})(t_i+t_{i+1}-2(t_{f}-t_0)) \\
    \frac{\sin((n-\Dot{\omega}_c)(t_{i+1}-t_0) + \Dot{\omega}_c(t_f-t_0) + u_{c,0})-\sin((n-\Dot{\omega}_c)(t_{i}-t_0) + \Dot{\omega}_c(t_f-t_0) + u_{c,0})}{n-\Dot{\omega}_c} \\
    \frac{\cos((n-\Dot{\omega}_c)(t_{i}-t_0) + \Dot{\omega}_c(t_f-t_0) + u_{c,0})-\cos((n-\Dot{\omega}_c)(t_{i+1}-t_0) + \Dot{\omega}_c(t_f-t_0) + u_{c,0})}{n-\Dot{\omega}_c} \\
    0 \\
    0.5\Phi_{\text{61}}^{\text{J2}}(t_i-t_{i+1})(t_i+t_{i+1}-2(t_{f}-t_0)) \\
    \end{bmatrix}
    \Delta B^{\rho}(t_i) \nonumber \\
    &= \boldsymbol{\Gamma}^{\text{drag}}(t_{i})  \Delta B^{\rho}(t_i) = \boldsymbol{\Gamma}^{\text{drag}}(t_{i})  \Delta B^{\rho}_i \nonumber
\end{align}
\noindent where the substitution for chief mean argument of latitude $u_c = n t + u_{c,0}$ is used and $u_{c,0} = \omega_{c,0}+M_{c,0}$ defines the initial mean argument of latitude. Here, the formulation has been rearranged the isolate the control variable, $\Delta B^{\rho}$, from the variables dependent on the spacecraft's state, which are concatenated into the auxiliary matrix $\boldsymbol{\Gamma}^{\text{drag}}$. In many cases, this integration is analytically intractable, due to time-varying parameters in multiple periodic terms. Two key assumptions remedy this issue. First, atmospheric density is assumed to be constant over each time step at its value at the beginning of the time step. This neither requires nor excludes the possibility of the density value being the same for the chief and deputy. Second, the differential ballistic coefficient is assumed to be held constant over the time step, such that the entire augmented differential ballistic coefficient is no longer time-dependent and reduces to $\Delta B^{\rho}_i$ for each time step. This effectively places a lower bound constraint on drag maneuver length, dictating the shortest drag maneuver that can be conducted. 

Incorporating the affine LTV formulation of drag control into Eq. \ref{eq:general expanded control problem} produces the following formulation
\begin{align} \label{eq: hybrid control problem}
    \text{minimize  } &c = \sum_{j=1}^{k} ||\boldsymbol{u}^{\text{propulsive}}_{j}||_{2} \\
    \text{subject to  } &\Delta \delta \boldsymbol{\alpha} = \sum_{j=1}^{k} \boldsymbol{\Gamma}(t_{j}) \boldsymbol{u}^{\text{propulsive}}_{j} + \sum_{i=1}^{m} \boldsymbol{\Gamma}^{\text{drag}}(t_{i}) \Delta B^{\rho}_{i} \qquad t_i \in T \text{ , } t_j \in T \nonumber\\
    &\sum_{i=1}^{m} t_{i+1}-t_i = t_f-t_0 \nonumber \\
    &\Delta B_{i,\text{min}}^{\rho} \leq \Delta B_i^{\rho}(t) \leq \Delta B_{i,\text{max}}^{\rho} \qquad \forall \; t_i \nonumber
\end{align}
\noindent which reduces to a norm-like cost function and fully LTV dynamics, making the problem compatible with the reachable set theory method proposed in [\citen{chernick}]. This formulation neglects the rotational dynamics of both the chief and deputy spacecraft and as a result assumes changes in spacecraft attitude occur instantaneously. This must be taken into account during control design when discretizing the control window for drag maneuvering, as the minimum drag maneuver length cannot be smaller than the time required to change spacecraft differential attitude from its maximum to minimum value. The key benefit of adding differential drag control is that it incurs little or no propellant cost, and the cost of drag maneuvers can be removed from the cost function. The resulting optimal control problem is almost identical to the impulsive control problem posed in Eq. 2 in [\citen{chernick}] by Chernick with only the addition of drag control in the dynamics constraint. It is easy to see that if the effect of differential drag on the pseudostate is known, then it can be subtracted from the original pseudostate of the reconfiguration. Any pseudostate-based optimal control solver can solve this reduced problem without modification.

\section{Hybrid Reachable Sets}
In general, reachable set theory analyzes the maximum effect of a given maneuver on the spacecraft's pseudostate. In the impulsive case, the shape of the reachable set is formed by sampling all possible directions for a maneuver of arbitrary cost. The convex hull of these points establishes the minimum cost required for a convex combination of maneuver magnitudes to reach a desired pseudostate. Because differential drag is a zero-cost extended maneuver and has control magnitude limits, the optimal policy for differential drag is bang-bang control, assuming the desired final state is not reachable with drag-control alone. The differential drag reachable set can be visualized by varying the length of a single drag maneuver within the desired control window. The goal of the reachable set theory analysis is then to find the optimal drag maneuver length that minimizes the remaining cost of the reconfiguration control problem. 

\begin{figure} [h]
    \centering
    \subfloat{\includegraphics[width=0.43\textwidth]{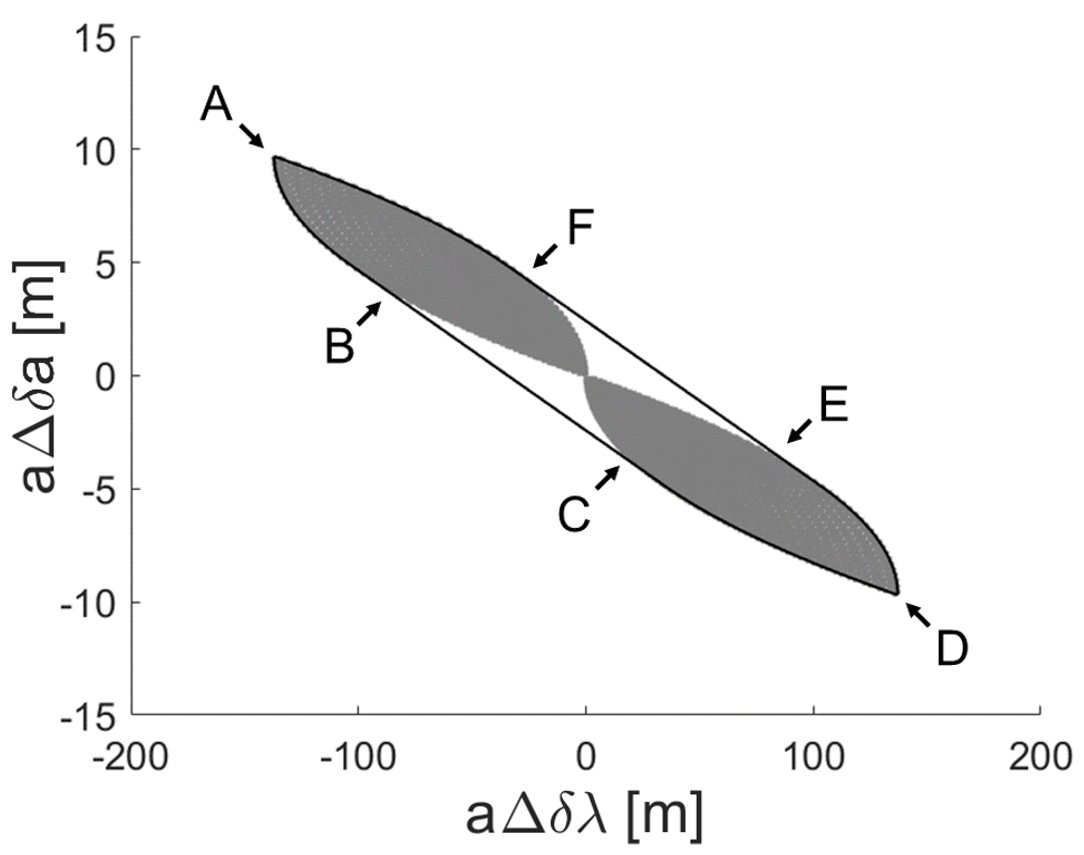}}
    \subfloat{\includegraphics[width=0.43\textwidth]{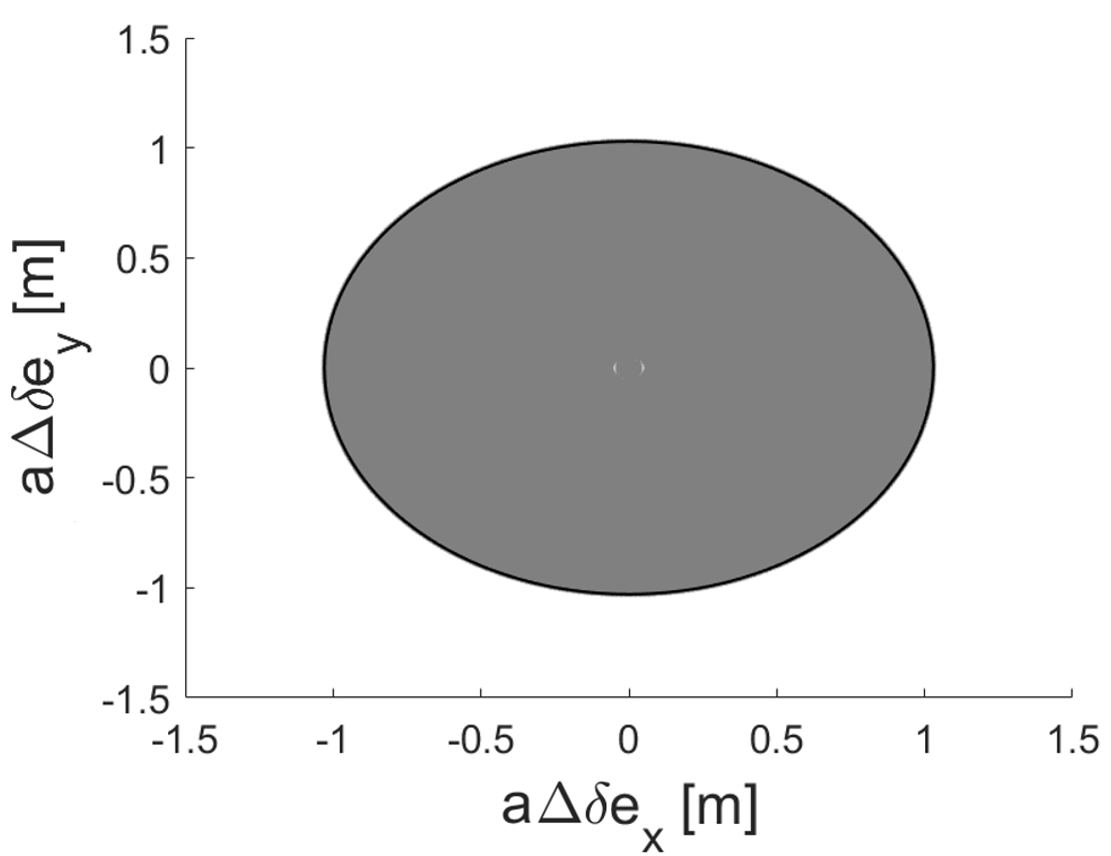}}
    \caption{Differential drag reachable sets under constant atmospheric density in the $\Delta \delta \boldsymbol{a}$ plane (left) and the $\Delta \delta \boldsymbol{e}$ plane (right).}
    \label{fig:RS drag}
    \vspace{-3mm}
\end{figure}

While differential drag affects the out-of-plane pseudostate in Eq. \ref{eq: integrated drag}, the effect on the $\delta i_y$ pseudostate is several magnitudes smaller than any of the in-plane ROE pseudostates and can be safely neglected. Therefore, only the in-plane ROE need to be analyzed. The shape of the differential drag reachable sets are found in Figure \ref{fig:RS drag} through Eq. \ref{eq: integrated drag}, by discretizing a three orbit control window and sampling all viable combinations of $t_i$  and $t_{i+1}$ ($t_{i+1} \geq t_i$) for both the positive and negative limits of $\Delta B^{\rho}$.

The visualization assumes equal positive and negative limits on $\Delta B^{\rho}$ and constant atmospheric density across the formation and throughout the control window. It is worth noting that the constant density assumption is required to maintain the consistency of the differential drag reachable set such that closed-form solutions can be found. Given that drag is posed as an affine approximation of many zero-cost impulsive maneuvers, the convex hull of these sample points is the reachable set for differential drag, indicated by the black lines in Figure \ref{fig:RS drag}. Figure \ref{fig:RS drag} also indicates important points of the reachable set in the $\Delta \delta \boldsymbol{a}$ plane for the following hybrid reachable set analysis.

\begin{figure} [h]
    \centering
    \subfloat{\includegraphics[width=0.43\textwidth]{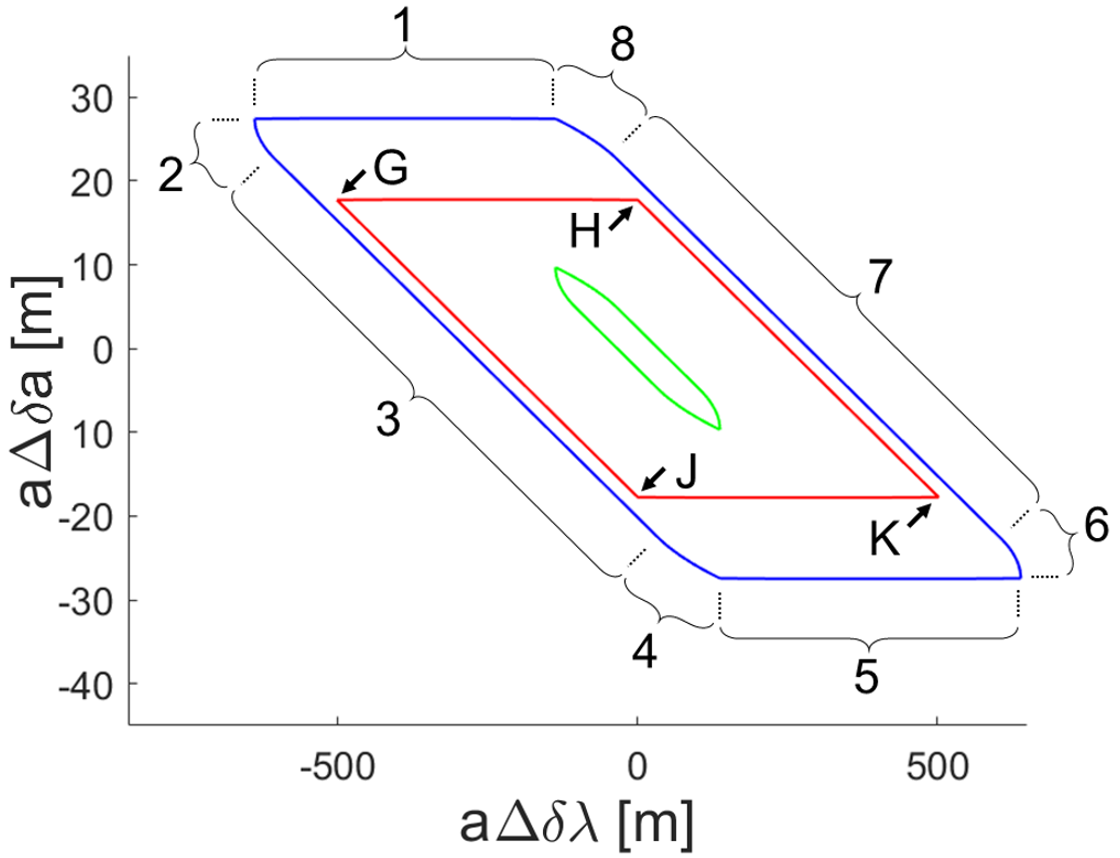}}
    \subfloat{\includegraphics[width=0.43\textwidth]{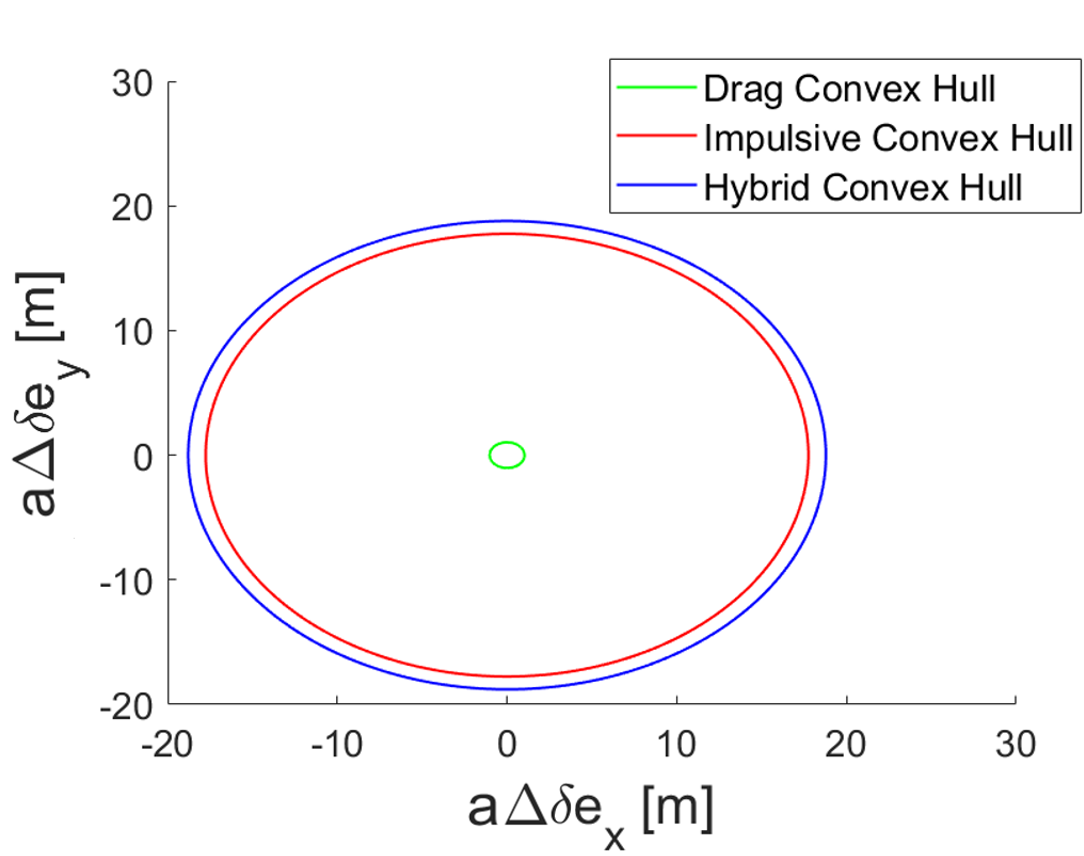}}
    \caption{Comparison of the reachable sets for different types of control under constant atmospheric density in the $\Delta \delta \boldsymbol{a}$ plane (left) and the $\Delta \delta \boldsymbol{e}$ plane (right).}
    \label{fig:RS hybrid}
    \vspace{-3mm}
\end{figure}

Hybrid control is analyzed by taking all possible combinations of these sample points with the near-circular impulsive reachable sets found by Chernick. The impulsive reachable sets are visualized by sampling impulsive maneuvers of magnitude 0.01m/s in all RTN directions and plugging these maneuvers into Eq. \ref{eq: RST single maneuver}. The convex hull of the resulting pseudostates is the impulsive reachable set. The hybrid reachable set is found by taking all combinations of sampled points between both sets and summing their effect to produce hybrid samples. The convex hull of these hybrid samples is the hybrid reachable set. Figure \ref{fig:RS hybrid} compares the varying shapes of the impulsive, differential drag, and hybrid reachable sets with important relevant contours and points indicated.
\vspace{-3mm}

The simplicity of the shape of the impulsive reachable set allows some important limits to be placed on the range of possible drag profiles, also known as the differential drag maneuver plan or the modulation of the differential ballistic coefficient over the control window. Portions of the hybrid reachable set contour in the $\Delta \delta \boldsymbol{a}$ plane are clearly dictated by either the shape of the impulsive or differential drag reachable set. Contours 1, 3, 5, and 7 are all straight linear lines, like the impulsive set, while Contours 2, 4, 6, and 8 are all nonlinear curves with the same shape as corresponding sections of the drag set. The hybrid contour in the $\Delta \delta \boldsymbol{e}$ plane is a circle, the same as both the impulsive and drag reachable sets. The critical observation to make about the hybrid reachable set is that the reconfiguration dominance cases for the impulsive-only problem ($\delta a$, $\delta \lambda$, $\delta e$) are the same as the hybrid problem, as drag maneuvers contribute no cost. This is best illustrated graphically. Consider Contour 2 indicated in Figure \ref{fig:RS hybrid}. This contour is created by summing the pseudostates in Segment A-B in the drag reachable set with that at Point G on the impulsive reachable set. When considering the combined cost, Segment A-B contributes no cost since it is composed of drag maneuvers, and Point G corresponds to the exact point on the impulsive reachable set where the $\delta a$ and $\delta \lambda$ dominance cases are equal. Therefore, a pseudostate that intersects Contour 2 indicates equal delta-v minimums for the $\delta a$ and $\delta \lambda$ dominance cases. All drag maneuver lengths in the contour are equally optimal, since the cost of the drag profile is zero regardless of the length of the drag maneuver. In contrast, consider Contour 1. This contour is created by summing the pseudostate at Point A in the drag reachable set with those in Segment G-H on the impulsive reachable set. Again considering cost, Point A contributes no cost as a drag maneuver, and Segment G-H encompasses the complete $\delta a$ dominance case for $\Delta \delta a > 0$. Therefore, Contour 1 only corresponds with one impulsive dominance case, $\delta a$, as well as one optimal drag profile length at Point A. Additionally, recall that this maneuver length at Point A was optimal for all of Contour 2. By extending this thought process to the full hybrid reachable set, the optimal drag maneuver lengths that correspond to contours with only one optimal drag maneuver length and one impulsive dominance case also achieve optimal cost in the regions with a range of optimal drag maneuver lengths and multiple equal impulsive dominance cases. Therefore, the only drag profiles required to reach global optimality are those that correspond to a single impulsive dominance case. The optimal drag control profiles can be analytically derived for each one of these points, such that all possible optimal control profile variations are known a priori. This is done in the following section.

\section{Optimal Drag Profiles}
\begin{figure} [htb]
    \centering
    \subfloat{\includegraphics[width=0.5\textwidth]{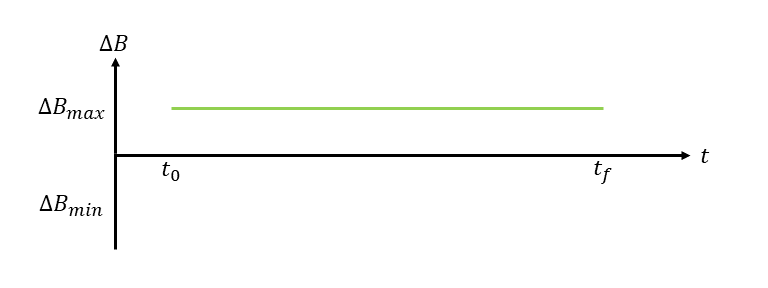}}
    \subfloat{\includegraphics[width=0.5\textwidth]{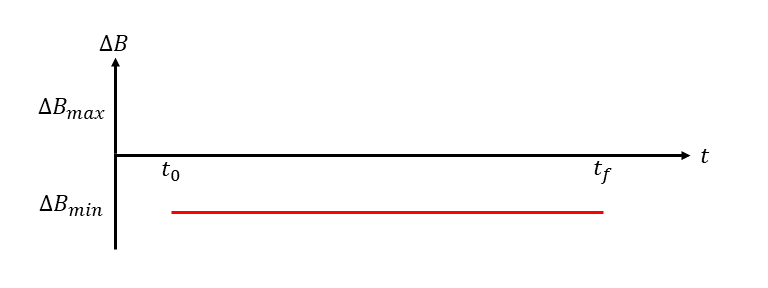}}
    \caption{Optimal drag profiles for $\delta a$ dominant reconfigurations given $\Delta \delta a > 0$ (left) or $\Delta \delta a < 0$ (right). Green drag profile corresponds to positive $\Delta B$, and red drag profile corresponds to negative $\Delta B$.}
    \label{fig:da drag profiles}
\end{figure}

As mentioned previously, only three dominance cases exist for in-plane ROE in near-circular orbit: $\delta a$, $\delta\lambda$, and $\delta e$. The strictly $\delta a$ dominant regions (Contours 1 and 5, cf. Fig. \ref{fig:RS hybrid}) correspond to a single drag maneuver at Points A and D (cf. Fig. \ref{fig:RS drag}). The optimal drag maneuver spans the entire control window in the direction of the desired $\delta a$ pseudostate, and both options are displayed in Figure \ref{fig:da drag profiles}.
\vspace{-4mm}

\begin{figure} [htb]
    \centering
    \subfloat{\includegraphics[width=0.5\textwidth]{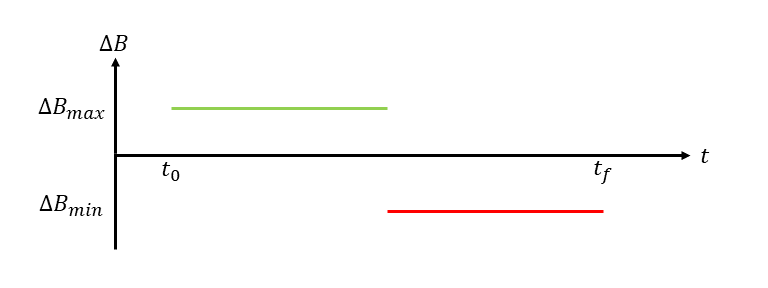}}
    \subfloat{\includegraphics[width=0.5\textwidth]{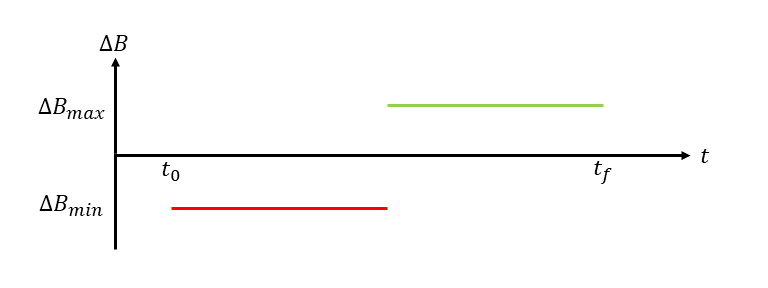}}
    \caption{Optimal drag profiles for $\delta \lambda$ dominant reconfigurations given $\Delta \delta \lambda < 0$ (left) or $\Delta \delta \lambda > 0$ (right). Green drag profile corresponds to positive $\Delta B$, and red drag profile corresponds to negative $\Delta B$.}
    \label{fig:dl drag profiles}
\end{figure}

For the strictly $\delta\lambda$ dominant regions (Contours 3 and 7, cf. Fig. \ref{fig:RS hybrid}), the corresponding drag maneuver occurs at Points B, C, E, and F (cf. Fig. \ref{fig:RS drag}) where the convex hull of the drag reachable set is tangent to the $\delta\lambda$ dominant regions of the impulsive set. The slope of the $\delta\lambda$ dominant regions was found in Table 12 in [\citen{chernick}] by Chernick, defined as
\begin{equation} \label{eq: dl slope}
    m = \frac{2}{\Hat{\Phi}_{\text{21}}^{\text{J2}}(t_f-t_0)}
\end{equation}
\noindent where $m$ is the slope and $\Hat{\Phi}_{\text{21}}^{\text{J2}} = \Phi_{\text{21}}^{\text{J2}}/(t_f-t)$ is the time independent component of $\Phi_{\text{21}}^{\text{J2}}$. The parameterized slope of the drag convex hull in the $\Delta \delta \boldsymbol{a}$ plane is defined from Eq. \ref{eq: integrated drag} as
\begin{equation}
    \frac{\Delta\delta a}{\Delta\delta \lambda} = \frac{\Delta t}{0.5\Hat{\Phi}_{\text{21}}^{\text{J2}}(-\Delta t)(2t+\Delta t-2(t_{f}-t_0))} \approx \frac{2}{\Hat{\Phi}_{\text{21}}^{\text{J2}}(2(t_{f}-t_0)-2t)} 
\end{equation}
\noindent where $t_i=t$ and $t$ is the maneuver timing, $t_{i+1}=t+\Delta t$ and $\Delta t$ is the length of an infinitesimally small drag maneuver, and $\Delta \delta a$ and $\Delta \delta \lambda$ are changes in pseudostate due to drag. Setting this slope equal to the slope of the $\delta\lambda$ dominant region solves the tangent point between the drag convex hull and the $\delta \lambda$ dominant region as
\begin{align}
    \frac{2}{\Hat{\Phi}_{\text{21}}^{\text{J2}}(2(t_{f}-t_0)-2t)} = m =& \frac{2}{\Hat{\Phi}_{\text{21}}^{\text{J2}}(t_f-t_0)} \\
    t = \frac{1}{2}(t_f-t_0)& \nonumber
\end{align}
\noindent The nonlinear contours of the drag convex hull in the $\Delta \delta \boldsymbol{a}$ plane are shaped by maneuvers that either start at $t_0$ or end at $t_f$. Two optimal maneuvers result: a maneuver in the direction of the desired $\delta \lambda$ pseudostate that spans the first half of the control window and a maneuver in the opposite direction of the desired $\delta \lambda$ pseudostate that spans the second half of the control window. Both maneuvers are achievable in the same control window, and the optimal drag profiles are shown in Figure \ref{fig:dl drag profiles}.

\begin{figure} [htb]
    \centering
    \subfloat{\includegraphics[width=0.5\textwidth]{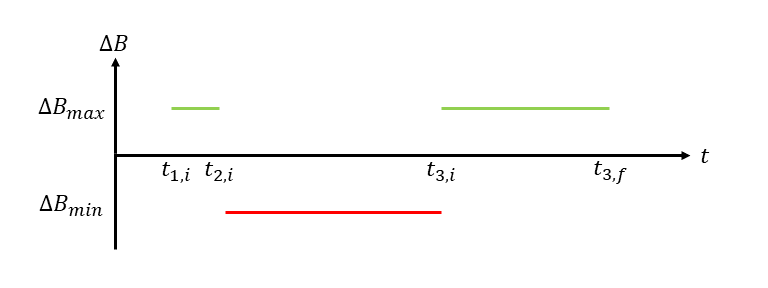}}
    \subfloat{\includegraphics[width=0.5\textwidth]{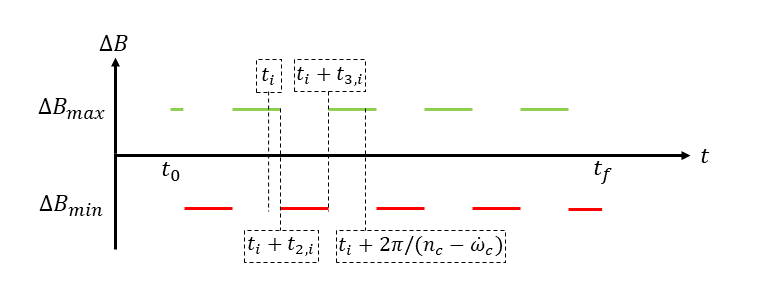}}
    \caption{Optimal drag profile for $\delta e$ dominant reconfigurations for a single cycle (left) and the full control window with a single cycle indicated (right). Green drag profile corresponds to positive $\Delta B$, and red drag profile corresponds to negative $\Delta B$.}
    \label{fig:de drag profiles}
\end{figure}

Similar to propulsive maneuvers, the effect of differential drag on the $\Delta \delta \boldsymbol{e}$ plane is periodic and circular. The drag reachable set shape matches the impulsive set shape in all directions, and while this means that the optimal profile varies over the entire hybrid set, several basic characteristics of the profile remain constant. First, the magnitude of the differential drag $\Delta \delta \boldsymbol{e}$ vector is analyzed by considering an arbitrary drag maneuver beginning at time $t$ until time $t+\Delta t$, given again from Eq. \ref{eq: integrated drag} as
\begin{align} 
    \Delta\delta e_x = \frac{a_c n_c \Delta B^{\rho}}{n_c-\Dot{\omega}_c}& [\sin((n-\Dot{\omega}_c)(t+ \Delta t -t_0) + \Dot{\omega}_c(t_f-t_0) + u_{c,0}) \\ &-\sin((n-\Dot{\omega}_c)(t -t_0) + \Dot{\omega}_c(t_f-t_0) + u_{c,0})] \nonumber \\
    = \frac{a_c n_c \Delta B^{\rho}}{n_c-\Dot{\omega}_c}&[2\cos((n-\Dot{\omega}_c)(t+ \frac{1}{2}\Delta t -t_0) + \Dot{\omega}_c(t_f-t_0) + u_{c,0})\sin((n-\Dot{\omega}_c)\frac{1}{2}\Delta t)] \nonumber \\
    \Delta\delta e_y = \frac{a_c n_c \Delta B^{\rho}}{n_c-\Dot{\omega}_c}&[\cos((n-\Dot{\omega}_c)(t -t_0) + \Dot{\omega}_c(t_f-t_0) + u_{c,0}) \\
    &-\cos((n-\Dot{\omega}_c)(t+ \Delta t -t_0) + \Dot{\omega}_c(t_f-t_0) + u_{c,0})] \nonumber \\
    = \frac{a_c n_c \Delta B^{\rho}}{n_c-\Dot{\omega}_c}&[2\sin((n-\Dot{\omega}_c)(t+ \frac{1}{2}\Delta t -t_0) + \Dot{\omega}_c(t_f-t_0) + u_{c,0})\sin((n-\Dot{\omega}_c)\frac{1}{2}\Delta t)] \nonumber \\
    &||\Delta \delta \boldsymbol{e}|| = \frac{a_c n_c \Delta B^{\rho}}{n_c-\Dot{\omega}_c} 2\sin((n-\Dot{\omega}_c)\frac{1}{2}\Delta t) \label{eq: de magnitude}
\end{align}
\noindent where $n_c-\Dot{\omega}_{c}$ is slightly less than the mean motion $n_c$, and the drag control period will be defined as $T_{\text{drag}}=\frac{2\pi}{n_c-\Dot{\omega}_{c}}$. The periodic function in Eq. \ref{eq: de magnitude} reaches a maximum value at $\Delta t = \frac{\pi}{n_c-\Dot{\omega}_c}$, so the optimal drag length is half of the drag control period. Both negative and positive maneuvers can reach this magnitude, such that the optimal policy involves oscillating maneuvers of this length repeated throughout the control window. The exact timing of these oscillations is such that both maneuvers will act directly towards the desired $\delta \boldsymbol{e}$ pseudostate. To capture all possible combinations over a given drag control period, a three-drag scheme beginning at an arbitrary time $t_i$ is considered in closed-form, defined as
\begin{equation}
    \frac{\Delta \delta e_{y,desired}}{\Delta \delta e_{x,desired}} = \frac{c_{1,i}-c_{1,f}-(c_{2,i}-c_{2,f})+c_{3,i}-c_{3,f}}{s_{1,f}-s_{1,i}-(s_{2,f}-s_{2,i})+s_{3,f}-s_{3,i}}
\end{equation}
\noindent where $\Delta \delta e_{y,desired}$ and $\Delta \delta e_{x,desired}$ are the desired pseudostates, $s_{m,j} = \sin((n-\Dot{\omega}_c)t_{i} + \Dot{\omega}_c(t_f-t_0) + u_{c,0})$, $c_{m,j} = \cos((n-\Dot{\omega}_c)t_{i} + \Dot{\omega}_c(t_f-t_0) + u_{c,0})$, subscript $m$ indicates the number of the drag maneuver, and subscript $j$ indicates the timing of the beginning $i$ or end $f$ of the drag maneuver. Given that the cycle spans one drag control period, several simplifying substitutions can be made, given as
\begin{align}
    &t_{1,f} = t_{2,i} \text{ , } t_{2,f} = t_{3,i} \\
    &s_{1,i} = \sin(\Dot{\omega}_c(t_f-t_0) + u_{c,0}) = s_{3,f} = \sin(2\pi + \Dot{\omega}_c(t_f-t_0) + u_{c,0}) \\
    &c_{1,i} = \cos(\Dot{\omega}_c(t_f-t_0) + u_{c,0}) = c_{3,f} = \cos(2\pi + \Dot{\omega}_c(t_f-t_0) + u_{c,0}) \\
    &t_{3,i} = t_{2,i} + \frac{\pi}{n-\Dot{\omega}_c} \label{eq: de direction 1}
\end{align}
by assuming the maneuvers are bang-bang. The second or “middle” maneuver is set to the optimal drag length in Eq. \ref{eq: de magnitude} and to the opposite direction of the other maneuvers. Because this optimal length is half of the drag control period, the other two maneuvers are also optimal, as they overlap over the end of one cycle to the next. An example of this three-drag scheme cycle is visualized in Figure \ref{fig:de drag profiles}. The drag oscillation timings are finally defined as
\begin{align} 
    &\frac{\Delta \delta e_{y,desired}}{\Delta \delta e_{x,desired}} = \frac{-2c_{2,i}}{2s_{2,i}} = -\cot((n-\Dot{\omega}_c)(t_{2,i} -t_0) + \Dot{\omega}_c(t_f-t_0) + u_{c,0}) \\
    &t_{2,i}  = \frac{1}{(n-\Dot{\omega}_c)} \left( \cot^{-1} \left( -\frac{\Delta \delta e_{y,desired}}{\Delta \delta e_{x,desired}} \right) - \Dot{\omega}_c(t_f-t_0) + u_{c,0}) \right) +t_0 \label{eq: de direction 2} \\
    &t_{3,i} = t_{2,i} + \frac{\pi}{n-\Dot{\omega}_c} +t_0
\end{align}
\noindent This defines a repeating cycle over the entire control window. Since the beginning time of these oscillations is arbitrary, the first cycle begins immediately at the start of the control window. An example of the full optimal profile for $\delta e$ dominant reconfigurations is displayed in Figure \ref{fig:de drag profiles}. The exact timing of the oscillations in the optimal profile will vary depending on the direction of the $\Delta \delta \boldsymbol{e}_{\text{desired}}$ vector.

\section{Closed-form Solutions for Drag Control}
The previous section covers all optimal drag control profiles for reconfigurations with a single hybrid dominance case, meaning that the subtracting the pseudostate reached by the dominant drag profile from the original pseudostate does not change the problem's dominance case. Optimally solving the hybrid control problem is not as simple as applying just one dominant profile. For example, recalling that the delta-v minimum for each dominance case is the cost required to achieve that particular ROE's pseudostate, consider a reconfiguration problem with a slightly higher $\delta a$ delta-v minimum than $\delta \lambda$ delta-v minimum for propulsive-only control. The highest delta-v minimum drives the minimum cost to accomplish the full reconfiguration, so applying the entire $\delta a$ dominant drag profile may result in a reduced propulsive reconfiguration problem that is $\delta \lambda$ dominant. This implies that part of the $\delta \lambda$ dominant profile should have been used to decrease the $\delta \lambda$ delta-v minimum to match the $\delta a$ delta-v minimum, indicating the full $\delta a$ dominant drag profile is non-optimal with respect to the complete hybrid control problem. Furthermore, each drag profile only optimally minimizes the ROE of its corresponding dominance case. Pseudostates that are fully reachable with only drag control may require contributions from all three dominant profiles to accomplish the full in-plane reconfiguration at zero cost. To improve drag profile robustness and maintain optimality while bridging the gap between hybrid and drag-only maneuver planning, a flexible closed-form drag profile structure is constructed, given knowledge of the state space and the limited number of dominance cases. To distinguish this profile from the optimal profiles for a single hybrid dominance case, the new profile will be referred to as the composite drag profile. The optimal profiles from the previous section will be referred to collectively as the optimal single hybrid dominance case drag profiles and individually by their specific dominance case ($\delta a$, $\delta \lambda$, $\delta e$).

The general methodology of the composite drag profile is to find the combination of the three single hybrid dominance case drag profiles that result in the lowest overall delta-v minimum. The equations for the propulsive delta-v minima of each in-plane dominance case in near-circular orbit were derived by Chernick in [\citen{chernick}] as 
\begin{align} 
    \delta v_{\text{min}, \delta a} &= \frac{n_c}{2}a_c\Delta \delta a \label{eq: Chernick dv min da}\\
    \delta v_{\text{min}, \delta \lambda} &= \frac{\pm(m a_c\Delta \delta\lambda - a_c\Delta  \delta a)}{a_c\Delta  \delta a_0} \label{eq: Chernick dv min dl}\\
    \delta v_{\text{min}, \delta e} &= \frac{n_c}{2}||a_c \Delta \delta \boldsymbol{e}||_2 
    \label{eq: Chernick dv min de}
\end{align}
\noindent where $\Delta \delta a_0$ defines the pseudostate achieved by a unitary tangential impulse at $t_0$ and $\delta v$ denotes a propulsive delta-v minimum. These $\delta v_{\text{min}}$ change for hybrid control, as the propulsive cost depends on the pseudostate of the original problem subtracted by that reached by the drag control profile. Finding the optimal cost of the composite drag profile requires comparing the effect of one optimal single hybrid dominance case drag profile against the others. First, the effect of using a portion of the $\delta a$ dominant profile on $\delta v_{\text{min}, \delta a}$ is found through Eq. \ref{eq: integrated drag} as
\begin{align} \label{eq: drag da ddv}
    \Delta \delta v_{\delta a} &= \frac{n_c}{2}a_c\Delta \delta a_{\text{drag}} = \frac{1}{2}a_{c}^{2}n_{c}^{2} \Delta B_{\text{approx}}^{\rho} \Delta t_{\delta a}
\end{align}
\noindent where $\Delta B_{\text{approx}}^{\rho} = \frac{1}{2}(||\Delta B_{\text{max}}^{\rho}|| + ||\Delta B_{\text{min}}^{\rho}||)$ approximates the augmented differential ballistic coefficient as its mean maximum magnitude, $\Delta \delta v$ denotes the change in the delta-v minimum due to using a drag profile, $\Delta \delta a_{\text{drag}}$ defines the change in pseudostate due to using the $\delta a$ dominant drag profile, and $\Delta t_{\delta a}$ is the length of the $\delta a$ dominant drag profile used. It is important to note that $\Delta \delta a_{\text{drag}}$ does not depend on when the $\delta a$ dominant profile is used in the control window, only the length of the dominant profile that is used. Next, the effect of using a portion of the $\delta \lambda$ dominant profile on $\delta v_{\text{min}, \delta \lambda}$ is found through Eq. \ref{eq: integrated drag} as
\begin{align} \label{eq: drag dl ddv}
    \Delta \delta v_{\delta \lambda} &= \frac{\pm(m a_c\Delta \delta\lambda_{\text{drag}} - a_c\Delta  \delta a_{\text{drag}})}{a_c\Delta  \delta a_0} \\
    &= \pm a_{c}^{2}n_c \Delta B^{\rho}\frac{\frac{1}{2}\Hat{\Phi}_{21}^{J_2}m(2 (t_f-t_0) \Delta t - 2t\Delta t - \Delta t^2) - \Delta t}{a_c \Delta \delta a_0} \nonumber \\
    &= 2a_{c}^{2}n_c \Delta B_{\text{approx}}^{\rho}\frac{-\frac{1}{2}\Hat{\Phi}_{21}^{J_2}m(\Delta t_{\delta \lambda}^2) + \Delta t_{\delta \lambda}}{a_c \Delta \delta a_0} \nonumber
\end{align}
where $t$ and $\Delta t$ are the start time and length of any arbitrary drag maneuver in the control window and $\Delta t_{\delta \lambda}$ is the length of the $\delta \lambda$ dominant profile used at the start and end of the control window. The $\delta \lambda$ dominant profile is used optimally at the start and end of the control window where any given drag profile has its largest effect on $\delta v_{\text{min}, \delta \lambda}$. Finally, effect of using a portion of the $\delta e$ dominant profile on $\delta v_{\text{min}, \delta e}$ is found through Eq. \ref{eq: integrated drag} as
\begin{align} \label{eq: drag de ddv}
    \Delta \delta v_{\delta e} &= \frac{n_c}{2}||a_c \Delta \delta \boldsymbol{e}_{\text{drag}}||_2 \approx \frac{n_c}{2} \left(\frac{2a_c n_c \Delta B^{\rho}}{n_c-\Dot{\omega}_c} - 0 \right) \frac{n_c-\Dot{\omega}_c}{\pi} \Delta t_{\delta e} = \frac{a_{c}^{2}n_{c}^{2}\Delta B_{\text{approx}}^{\rho}}{\pi} \Delta t_{\delta e} 
\end{align}
where $\Delta t_{\delta e}$ is the length of the $\delta e$ dominant profile used. The periodic evolution of $\Delta \delta \boldsymbol{e}_{\text{drag}}$ is approximated as linear using the maximum magnitude per half drag control period found in Eq. \ref{eq: de magnitude} to simplify periodic functions to an analytically tractable linear model. It is again important to note that $\Delta \delta \boldsymbol{e}_{\text{drag}}$ does not depend on when the $\delta e$ dominant profile is used in the control window, only the length of the dominant profile that is used.

To compare the effect of using portions of one single hybrid dominance case drag profile against the others, the composite drag profile initializes as the $\delta a$ dominant profile, and the impact of replacing a portion of this profile with the $\delta \lambda$ or $\delta e$ dominant profile is evaluated analytically. This choice of initialization is advantageous, as the lack of oscillating maneuvers will cancel out any effect on $\delta v_{\text{min},\delta e}$ and drag profiles arranged symmetrically about the midpoint of the control window have a net-zero effect on $\delta v_{\text{min}, \delta \lambda}$. This effectively decouples the effect of each single hybrid dominance case drag profile and corresponding delta-v minimum from the other two, allowing direct comparisons between just two delta-v minima at a time. The equation for $\delta v_{\text{min},\delta a}$ is updated from Eq. \ref{eq: Chernick dv min da} for the initialization, given as
\begin{equation} \label{eq: updated dv min da}
    \delta v_{\text{min},\delta a} = \frac{n_c}{2}a_c\Delta \delta a - \frac{1}{2}a_{c}^{2}n_{c}^{2} \Delta B_{\text{approx}}^{\rho}(t_f-t_0)
\end{equation}
\noindent Replacing the $\delta a$ dominant profile with the $\delta \lambda$ dominant profile such that $\delta v_{\text{min}, \delta a}=\delta v_{\text{min}, \delta \lambda}$ is formulated through Eqs. \ref{eq: drag dl ddv} and \ref{eq: drag da ddv} as
\begin{align} 
    \delta v_{\text{min},\delta \lambda} -& \Delta \delta v_{\delta \lambda} = \delta v_{\text{min},\delta a} + 2\Delta \delta v_{\delta a} \\
    \delta v_{\text{min},\delta \lambda} - 2a_{c}^{2}n_c \Delta B_{\text{approx}}^{\rho}&\frac{-\frac{1}{2}\Hat{\Phi}_{21}^{J_2}m \Delta t_{\delta \lambda/\delta a}^2 + \Delta t_{\delta \lambda/\delta a}}{a_c \Delta \delta a_0} \\
    &= \delta v_{\text{min},\delta a} + 2\left( \frac{1}{2}a_{c}^{2}n_{c}^{2} \Delta B_{\text{approx}}^{\rho} \Delta t_{\delta \lambda/\delta a} \right) \nonumber \\
    \frac{-\Hat{\Phi}_{21}^{J_2}ma_{c}^{2}n_c \Delta B_{\text{approx}}^{\rho}}{\Delta a_c \delta a_0}\Delta t_{\delta \lambda/\delta a}^2 &+ 
    \left(2\frac{a_{c}^{2}n_c \Delta B_{\text{approx}}^{\rho}}{\Delta a_c \delta a_0} + 
    a_{c}^{2}n_{c}^{2} \Delta B_{\text{approx}}^{\rho} \right) \Delta t_{\delta \lambda/\delta a}  \label{eq: dl,da matching} \\
    & + \delta v_{\text{min},\delta a} - \delta v_{\text{min},\delta \lambda} = 0 \nonumber
\end{align}
\noindent where $\Delta t_{\delta \lambda/\delta a}$ is the length of the $\delta \lambda$ dominant profile that will replace the $\delta a$ dominant profile at both the start and end of the control window. The relationship between $\delta v_{\text{min},\delta a}$ and $\delta v_{\text{min},\delta \lambda}$ is quadratic, requiring root-finding to find if a valid solution exists. If root-finding does not result in any real or positive roots, then no change should be made to the profile and $\Delta t_{\delta \lambda/\delta a} = 0$. To prevent the $\delta \lambda$ dominant profile from overshooting the $\delta \lambda$ pseudostate, the portion of the $\delta \lambda$ dominant profile required for $\delta v_{\text{min}, \delta \lambda}=0$ is found as
\begin{equation}
    \delta v_{\text{min},\delta \lambda} - \Delta \delta v_{\delta \lambda} = 0 \rightarrow \frac{-\Hat{\Phi}_{21}^{J_2}ma_{c}^{2}n_c \Delta B_{\text{approx}}^{\rho}}{\Delta a_c \delta a_0}\Delta t_{\delta \lambda/0}^2 + 
    \left(2\frac{a_{c}^{2}n_c \Delta B_{\text{approx}}^{\rho}}{a_c \Delta \delta a_0} \right) \Delta t_{\delta \lambda/0} - \delta v_{\text{min},\delta \lambda} = 0 \label{eq: dl,0 matching}
\end{equation}
\noindent where $\Delta t_{\delta \lambda/0}$ is the maximum length of the $\delta \lambda$ dominant profile that can replace the $\delta a$ dominant profile at both the start and end of the control window.

Next, replacing the $\delta a$ dominant profile with the $\delta e$ dominant profile such that $\delta v_{\text{min}, \delta a}=\delta v_{\text{min}, \delta e}$ is formulated through Eqs. \ref{eq: drag de ddv} and \ref{eq: drag da ddv} as
\begin{align} 
    \delta v_{\text{min},\delta e} - \Delta \delta v_{\delta e} = & \delta  v_{\text{min},\delta a} + 2\Delta \delta v_{\delta a} \\
    \delta v_{\text{min},\delta e} - \frac{a_{c}^{2}n_{c}^{2}\Delta B_{\text{approx}}^{\rho}}{\pi} \Delta t_{\delta e/\delta a} = &
    \delta v_{\text{min},\delta a} + 2\left( \frac{1}{2}a_{c}^{2}n_{c}^{2} \Delta B_{\text{approx}}^{\rho} \left(\frac{1}{2}\Delta t_{\delta e/\delta a} \right) \right) \\
    \Delta t_{\delta e/\delta a} = &\frac{\delta v_{\text{min},\delta e} - \delta v_{\text{min},\delta a}}{\frac{a_{c}^{2}n_{c}^{2}\Delta B_{\text{approx}}^{\rho}}{\pi} + \frac{1}{2}a_{c}^{2}n_{c}^{2} \Delta B_{\text{approx}}^{\rho}} \label{eq: de,da matching}
\end{align}
\noindent where $\Delta t_{\delta e/\delta a}$ is the length of the $\delta e$ dominant profile that will replace the $\delta a$ dominant profile both before and after the midpoint of the control window. If Eq. \ref{eq: de,da matching} results in a negative value of $\Delta t_{\delta e/\delta a}$, then using the $\delta e$ dominant profile will only raise the overall delta-v minimum, and $\Delta t_{\delta e/\delta a}=0$. The $\delta e$ dominant drag profile is centered on the middle of the control window because the $\delta \lambda$ dominant profile must be at the start and end of the control window to be optimal, and the $\delta e$ dominant profile does not have this time-based constraint. This placement maximizes the lengths of both the $\delta \lambda$ and $\delta e$ dominant profiles that can be used without overlap. To prevent the $\delta e$ dominant profile from overshooting the $\delta e$ pseudostate, the portion of the $\delta e$ dominant profile required for $\delta v_{\text{min}, \delta e}=0$ is found as
\begin{equation}
    \delta v_{\text{min},\delta e} - \Delta \delta v_{\delta e} = 0 \rightarrow \Delta t_{\delta e/0} = \frac{\delta v_{\text{min},\delta e} }{\frac{a_{c}^{2}n_{c}^{2}\Delta B_{\text{approx}}^{\rho}}{\pi}} \label{eq: de,0 matching}
\end{equation}
\noindent where $\Delta t_{\delta e/0}$ is the maximum length of the $\delta e$ dominant profile that can replace the $\delta a$ dominant profile both before and after the midpoint of the control window. 

If the sum of $\Delta t_{\delta \lambda/\delta a}$ and $\Delta t_{\delta e/\delta a}$ is greater than half of the control window, then the entire $\delta a$ dominant profile should be replaced by a combination of the $\delta \lambda$ and $\delta e$ dominant profiles. The effect of using the $\delta \lambda$ dominant profile is directly compared against the $\delta e$ dominant profile, formulated through Eqs. \ref{eq: drag de ddv} and \ref{eq: drag dl ddv} as
\begin{align} 
     \delta v_{\text{min},\delta e} - & \Delta \delta v_{\delta e} = \delta v_{\text{min},\delta \lambda} - \Delta \delta v_{\delta \lambda} \\
     \delta v_{\text{min},\delta e} - \frac{a_{c}^{2}n_{c}^{2}\Delta B_{\text{approx}}^{\rho}}{\pi} 
     2 &\left( \frac{1}{2}(t_f-t_0) - \Delta t_{\delta \lambda/\delta e} \right) = \delta v_{\text{min},\delta \lambda} \nonumber \\
     &- 2a_{c}^{2}n_c \Delta B_{\text{approx}}^{\rho} \frac{-\frac{1}{2}\Hat{\Phi}_{21}^{J_2}m(\Delta t_{\delta \lambda/\delta e}^2) + \Delta t_{\delta \lambda/\delta e}}{a_c \Delta \delta a_0} \\
     \frac{-\Hat{\Phi}_{21}^{J_2}ma_{c}^{2}n_c \Delta B_{\text{approx}}^{\rho}}{\Delta a_c \delta a_0}\Delta t_{\delta \lambda/\delta e}^2 + 
     &\left( 2\frac{a_{c}^{2}n_c \Delta B_{\text{approx}}^{\rho}}{\Delta a_c \delta a_0} + 
    2\frac{a_{c}^{2}n_{c}^{2}\Delta B_{\text{approx}}^{\rho}}{\pi} \right) \Delta t_{\delta \lambda/\delta e}  \label{eq: de,dl matching} \\
     &- \frac{a_{c}^{2}n_{c}^{2}\Delta B_{\text{approx}}^{\rho}}{\pi} (t_f-t_0) + \delta v_{\text{min},\delta e} - \delta v_{\text{min},\delta \lambda} = 0 \nonumber
\end{align}
\noindent where $\Delta t_{\delta \lambda/\delta e}$ is the length of the $\delta \lambda$ dominant profile at the start and end of the control window. The $\delta e$ dominant profile is used in the rest of the control window. Root-finding is again required. If there are no real roots, the full $\delta \lambda$ dominant profile should be used, and if no positive real roots are less than half of the control window, the full $\delta e$ dominant profile should be used.

If the sum of $\Delta t_{\delta \lambda/\delta a}$ and $\Delta t_{\delta e/\delta a}$ is not greater than half of the control window, then both replacements can be accomplished without overlap. All in-plane delta-v minima are updated for this new profile. The final step to form the full composite drag profile is to ensure that the remaining portion of the $\delta a$ dominant profile will not overshoot the $\delta a$ pseudostate. If $\delta v_{\text{min}, \delta a}<0$, then the length of the $\delta a$ dominant profile should be reduced, formulated through Eq. \ref{eq: drag da ddv} as
\begin{equation}
    \delta v_{\text{min},\delta a} - 2\Delta \delta v_{\delta a} = 0 \rightarrow \Delta t_{\delta a/0} = \frac{\delta v_{\text{min},\delta a}}{2\left( \frac{1}{2}a_{c}^{2}n_{c}^{2} \Delta B_{\text{approx}}^{\rho} \right)} \label{eq: da,0 matching}
\end{equation}
\noindent where $\Delta t_{\delta a/0}$ is the length of the $\delta a$ dominant profile that should be removed from the composite drag profile. This removal is done symmetrically about the midpoint of the control window to negate any effect on $\delta v_{\text{min}, \delta \lambda}$. A negative $\delta v_{\text{min}, \delta a}$ at this final step indicates that the full in-plane pseudostate is reachable with drag control alone, and this final accommodation bridges the gap between hybrid and drag-only control. 

Using the equations and process presented above, optimal or near-optimal drag profiles for both hybrid and drag-only control can be found in a single closed-form procedure. The full algorithm used to construct the composite drag profile is detailed in Algorithm \ref{alg: flex drag profile}, which encompasses Eqs. \ref{eq: drag da ddv}-\ref{eq: da,0 matching}. It is worth noting that the composite profile will result in the previously found single hybrid dominance case optimal drag profiles if one dominance case has a much higher delta-v minimum than the others. Once the drag profile is produced by Algorithm \ref{alg: flex drag profile}, the pseudostate reached by drag control can be calculated and subtracted from the initial pseudostate of the hybrid control problem. The propulsive control required to reach this reduced pseudostate can be calculated by any optimal propulsive control solver without any additional information from the drag profile. The most reasonable option to maintain computational efficiency and consistency of approach is the closed-form impulsive algorithm designed by Chernick and D'Amico. \cite{chernick}

\begin{algorithm} [h] 
    \caption{Composite Drag Profile Algorithm} \label{alg: flex drag profile}
    \KwIn{$\boldsymbol{\alpha}_c$ , $\Delta \delta \boldsymbol{\alpha}_{\text{desired}}$ , $\Delta B_{\text{max}}^{\rho}$ , $\Delta B_{\text{min}}^{\rho}$ , $t_0$ , $t_f$}
    \KwResult{$\boldsymbol{\Delta B}^{\rho}_{\text{profile}}$}
    Initialize $\boldsymbol{\Delta B}^{\rho}_{\text{profile}}$ as the optimal $\delta a$ drag profile \\
    Find optimal $\delta e$ drag profile switching times: $t_{2,i}$ , $t_{3,i}$ $\gets$ \text{Eqs.} \ref{eq: de direction 1} \text{ and } \ref{eq: de direction 2} \\
    $\delta v_{min,\delta a}$ , $\delta v_{min,\delta \lambda}$ , $\delta v_{min,\delta e}$ $\gets$ Eqs. \ref{eq: updated dv min da}, \ref{eq: Chernick dv min dl}, \ref{eq: Chernick dv min de} \\
    $\Delta t_{\delta \lambda} \gets $ max( min( Eq. \ref{eq: dl,da matching} real roots , Eq. \ref{eq: dl,0 matching} real roots ) , 0 ) \\
    $\Delta t_{\delta e} \gets $ max( min( \text{Eq.} \ref{eq: de,da matching} , \text{Eq.} \ref{eq: de,0 matching} ) , 0 ) \\
    \eIf{$\Delta t_{\delta \lambda} + \Delta t_{\delta e} > \frac{1}{2}(t_f-t_0)$}
    {$\Delta t_{\delta \lambda/\delta e} \gets $ \text{Eq.} \ref{eq: de,dl matching} roots\\
    \uIf{$\Delta t_{\delta \lambda/\delta e} > \frac{1}{2}(t_f-t_0)$ \text{ or no real roots exist}}
    {Apply full $\delta \lambda$ dominant profile}
    \uElseIf{$\Delta t_{\delta \lambda/\delta e} < 0$}
    {Apply full $\delta e$ dominant profile}
    \Else
    {Apply $\delta \lambda$ dominant profile until $\Delta t_{\delta \lambda/\delta e}$ after $t_0$ and before $t_f$ \\
    Apply $\delta e$ optimal profile in remaining control window}
    }
    {Apply $\Delta t_{\delta \lambda}$ of $\delta \lambda$ dominant profile after $t_0$ and before $t_f$ \\
    Apply $\Delta t_{\delta e}$ of $\delta e$ dominant profile before and after  $\frac{1}{2}(t_f-t_0)$ }
    $\Delta \delta \boldsymbol{\alpha}_{\text{post-drag}} = \Delta \delta \boldsymbol{\alpha}_{\text{desired}} - \sum_{i=1}^{m} \boldsymbol{\Gamma}^{\text{drag}}(t_{i})  \Delta B^{\rho}(t_i) dt$ \\
    Update $\delta v_{min,\delta a}$ , $\delta v_{min,\delta \lambda}$ , $\delta v_{min,\delta e}$ $\gets$ Eqs. \ref{eq: Chernick dv min da}, \ref{eq: Chernick dv min dl}, \ref{eq: Chernick dv min de} \\
    \If{$\delta v_{min,\delta a} < 0$}
    {Reduce $\delta a$ dominant drag profile by Eq. \ref{eq: da,0 matching} \\
    Update $\delta v_{min,\delta a}$ $\gets$ Eq. \ref{eq: Chernick dv min da}}
\end{algorithm}

\section{Validation and Simulations}

The effectiveness of the proposed approach is verified in two ways. Given that the original control problem in Eq. \ref{eq: hybrid control problem} is a linear program, generic interior point algorithms will provide baselines with which the closed-form solutions can be compared for delta-v optimality and consistency with expected drag profiles. Both the numerical and closed-form hybrid maneuver plans, as well as optimal propulsive-only plans, for a given desired reconfiguration are then simulated in a full-force orbital propagator that integrates the GVE (4th order Runge-Kutta) for perturbations including Earth's gravity (30x30 gravity model), atmospheric drag, solar radiation pressure, third body effects, and relativistic corrections. The orbital propagator takes as input the osculating orbital elements for each spacecraft, transformed from the mean orbital elements using the Brouwer transformation, and propagates each forward over a 10s time step. \cite{brouwer} The resulting osculating orbital elements are converted back into mean orbital elements for post-processing calculation of mean ROE. These simulations will provide a baseline comparison of final state accuracy for the proposed approach. The chief and deputy spacecraft used in these simulations are identical 3U CubeSats based on those that will be used on the SWARM-EX mission, a major motivator of this work. \cite{swarm-ex} These spacecraft are assumed to have a mass of 6kg, a maximum cross-sectional area of 0.09$\text{m}^2$, a minimum cross-sectional area of 0.01$\text{m}^2$, and a coefficient of drag of 1.5. The initial mean orbital elements of the chief spacecraft are provided in vector form as
\begin{equation}
    \boldsymbol{\alpha}_c = 
    \begin{bmatrix}
    a & e & i & \Omega & \omega & M
    \end{bmatrix}
     = 
     \begin{bmatrix}
     6798km & 0.001 & 51^{\circ} & 0^{\circ} & 0^{\circ} & 90^{\circ}
     \end{bmatrix}
\end{equation}


The atmospheric density model used is the NRLMSISE-00 model, which provides variable atmospheric density measurements based on Earth-Centered Earth-Fixed (ECEF) position, solar flux, and geomagnetic activity. \cite{nrlmsise00} To isolate the performance of the algorithm from the uncertainty of atmospheric density estimation, the pseudostate achieved by any given drag profile through Eq. \ref{eq: integrated drag} is calculated with perfect knowledge of the atmospheric density of each spacecraft, and corresponding augmented ballistic coefficient, at all simulation time steps. A 200s time step is used to discretize the control window for drag control calculations. The propulsive solver of the proposed closed-from hybrid control algorithm is the before-mentioned approach derived by Chernick et al., referred to as the Chernick solver. \cite{chernick} This will also provide the propulsive-only control solutions to gauge the cost savings of hybrid control over propulsive-only methods. Because the Chernick solutions have not been evaluated for extended control windows, the numerical solver created by Koenig et al., referred to as the Koenig solver, will also produce in-plane delta-v minima for comparison and optimality evaluations. \cite{impulse_koenig} Numerically-computed hybrid control solutions for the original control problem presented in Eq. \ref{eq: hybrid control problem} are provided by SDPT3 through CVX. \cite{boyd} Table \ref{tab: reconfigs} details the reconfigurations that will be used to evaluate the various control approaches over a 30 orbit control window. Because differential drag has a very small out-of-plane effect, only in-plane delta-v minima and trajectories are analyzed and plotted. However for completeness, corresponding out-of-plane final state errors are also reported.

\begin{table}[h!]
    \centering
    \begin{tabular}{r|c|c|c|c|c|c}
    \hline \hline
    Deputy ROE & $a_c\delta a [m]$ & $a_c\delta\lambda [m]$ & $a_c\delta e_x [m]$ & $a_c\delta e_y [m]$ & $a_c\delta i_x [m]$ & $a_c\delta i_y [m]$ \\
    \hline
    Case 1 Initial & -300 & -30000 & 250 & 1900 & 100 & 1800 \\
    \hline
    Case 2 Initial & 30 & -75000 & -20 & 1970 & 100 & 1800 \\
    \hline
    Case 3 Initial & 30 & -6000 & -20 & 1970 & 100 & 1800 \\
    \hline
    Case 4 Initial & -20 & -2000 & 250 & 1980 & 100 & 1800 \\
    \hline
    Desired & 0 & -5000 & 0 & 2000 & 0 & 2000 \\
    \hline \hline
    \end{tabular}
    \caption{Reconfigurations for control validation over a 30 orbit control window.}
    \label{tab: reconfigs}
    \vspace{-3mm}
\end{table}

\subsection{Drag Profile Optimality}

\begin{figure} [h!]
    \centering
    \subfloat{\includegraphics[width=0.47\textwidth]{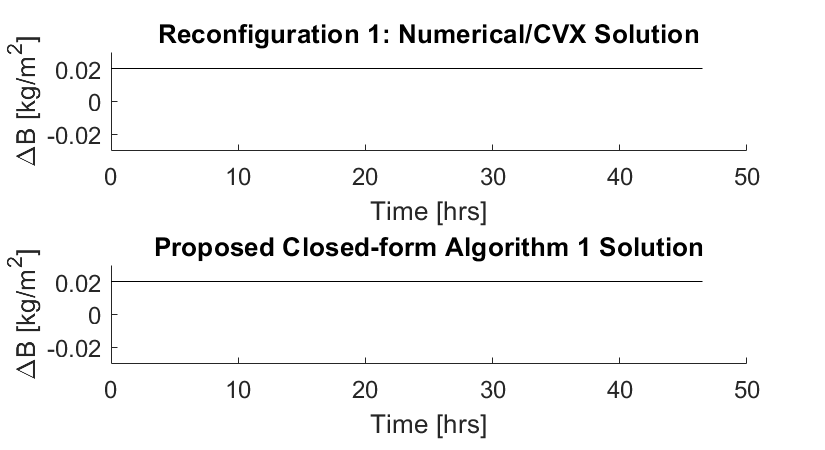}}
    \hfill
    \subfloat{\includegraphics[width=0.47\textwidth]{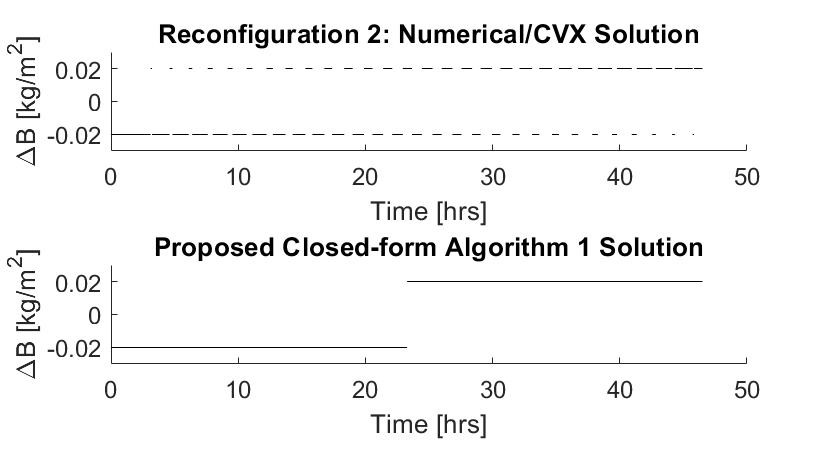}}\\
    \subfloat{\includegraphics[width=0.47\textwidth]{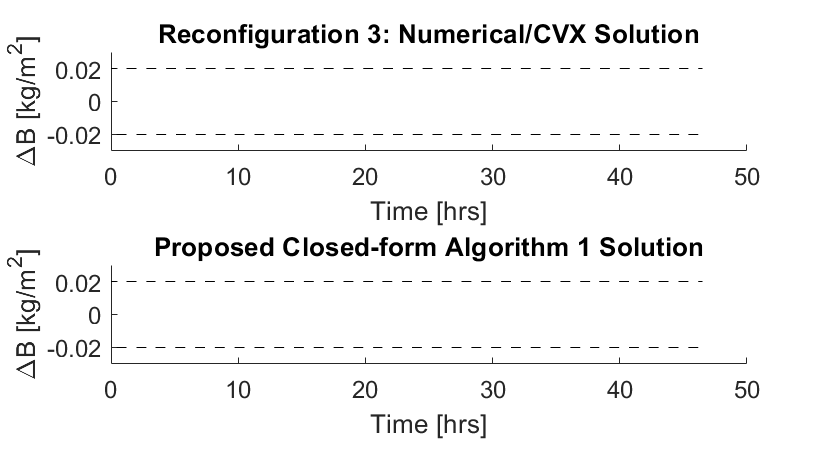}}
    \hfill
    \subfloat{\includegraphics[width=0.47\textwidth]{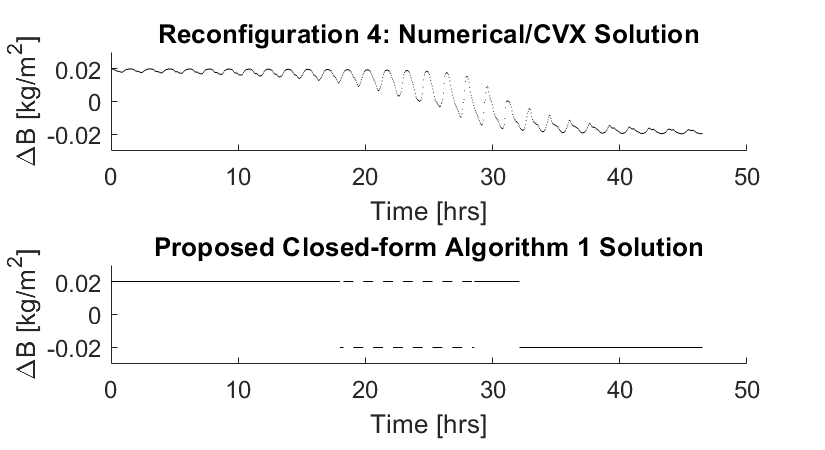}}
    \caption{Differential drag profiles for Table \ref{tab: reconfigs} validation reconfigurations.}
    \label{fig: drag profiles}
    \vspace{-5mm}
\end{figure}

The in-plane delta-v minima found by the differing control approaches for each reconfiguration in Table \ref{tab: reconfigs} are listed in Table \ref{tab: dvmin}, using  Eqs. \ref{eq: Chernick dv min da}-\ref{eq: Chernick dv min de} for the Chernick solver and the proposed closed-form hybrid approach. The differential drag control profiles calculated by the numerical solver and Algorithm \ref{alg: flex drag profile} are illustrated in Figure \ref{fig: drag profiles}.
\begin{table}[h!]
    \centering
    \begin{tabular}{c|c|c|c|c}
    \hline \hline
    Control Approach & Case 1 (m/s) & Case 2 (m/s) & Case 3 (m/s) & Case 4 (m/s) \\
    \hline
    Impulsive (Koenig) $\delta v_{\text{min}}$ & 0.1690 & 0.2988 & 0.1597 & 0.0229 \\
    \hline
    Impulsive (Chernick) $\delta v_{\text{min}}$ & 0.1690 & 0.2960 & 0.1597 & 0.0228 \\
    \hline
    Hybrid (Numerical/CVX) $\delta v_{\text{min}}$ & 0.1157 & 0.2712 & 0.1226 & 0.0000 \\
    \hline
    Hybrid (Closed-form) $\delta v_{\text{min}}$ & 0.1157 & 0.2693 & 0.1226 & 0.0036 \\
    \hline \hline
    \end{tabular}
    \caption{In-plane $\delta v_{min}$ found by various control approaches for Table \ref{tab: reconfigs} validation reconfigurations.}
    \label{tab: dvmin}
    \vspace{-3mm}
\end{table}

Reconfiguration 1 is $\delta a$ dominant with an in-plane hybrid $\delta v_{\text{min}}$ of 0.1157m/s found by Algorithm \ref{alg: flex drag profile} and Eq. \ref{eq: Chernick dv min da} with a drag profile that matches the $\delta a$ dominant drag profile. This hybrid solution provides 0.0398m/s delta-v cost savings (35.346\%) over closed-form propulsive-only methods. Both the closed-form hybrid and propulsive-only $\delta v_{\text{min}}$ match those found numerically with at least $10^{-4}$ precision. The numerical hybrid drag profile also results in exactly the same $\delta a$ dominant drag profile as the closed-form profile.

Reconfiguration 2 is $\delta \lambda$ dominant with an in-plane hybrid $\delta v_{\text{min}}$ of 0.2693m/s found by Algorithm \ref{alg: flex drag profile} and Eq. \ref{eq: Chernick dv min dl} with a drag profile that matches the $\delta \lambda$ dominant drag profile. The closed-form hybrid approach finds 0.0267m/s delta-v cost savings (9.020\%) over closed-form propulsive-only control. The Algorithm \ref{alg: flex drag profile} $\delta v_{\text{min}}$ results in 0.0019m/s error with respect to the hybrid numerical solution, which is less than the 0.0028m/s error between the closed-form and numerical propulsive-only solutions. Given that the Chernick closed-form impulsive solutions were proven to be globally optimal through comparison with the Koenig solver, the Algorithm \ref{alg: flex drag profile} solution is globally optimal. While the closed-form drag profile differs from the numerical solution for Reconfiguration 2, the delta-v minima in Table \ref{tab: dvmin} demonstrate that these profiles result in nearly equivalent propulsive costs. Furthermore, the profile found by Algorithm \ref{alg: flex drag profile} is significantly more realizable on a spacecraft when accounting for delays in attitude changes due to rotational dynamics.

Reconfiguration 3 is $\delta e$ dominant with an in-plane hybrid $\delta v_{\text{min}}$ of 0.1306m/s found by Algorithm \ref{alg: flex drag profile} and Eq. \ref{eq: Chernick dv min de} with a drag profile that matches the $\delta e$ dominant drag profile. The Algorithm \ref{alg: flex drag profile} solution produces 0.0291m/s delta-v cost savings (18.222\%) over closed-form propulsive-only control. Both the closed-form hybrid and propulsive-only $\delta v_{\text{min}}$ match those found numerically with at least $10^{-4}$ precision, and both drag profiles exactly match the $\delta e$ dominant drag profile.

Reconfiguration 4 corresponds to a nearly drag-only reconfiguration, requiring accommodations from all three single hybrid dominance case drag profiles. The closed-form hybrid solution results in a remaining $\delta e$ dominant problem with a 0.0036m/s in-plane $\delta v_{\text{min}}$ through Eq. \ref{eq: Chernick dv min de} after applying the drag profile. This is a significant difference from the approximately zero $\delta v_{\text{min}}$ found by the numerical solver ($\sim10^{-8}$m/s). However, this closed-form hybrid control cost is still greatly reduced from the propulsive-only solution, providing 0.0192m/s delta-v cost savings (84.211\%). The drag profiles found by the closed-form and numerical approaches also differ greatly. This is because the desired state is now reachable with drag-only control, such that bang-bang control is no longer a valid assumption. The bang-bang assumption is fundamental to the creation of closed-form drag profile solutions, and the only real improvements that could be made to Algorithm \ref{alg: flex drag profile} would be to use a more accurate approximation of $\Delta \delta v_{\delta e}$ in Eq. \ref{eq: drag de ddv} than the linear model used in this paper. This would increase the complexity of the closed-form algorithm for an unknown $\delta v_{\text{min}}$ benefit. Additionally, like the profiles from Reconfiguration 2, the profile found in closed-form is significantly more realizable on a spacecraft than the numerical solution for the drag profile.

\subsection{Hybrid Control Performance}

\begin{figure}[h!]
    \centering
    \subfloat{\includegraphics[width=0.44\textwidth]{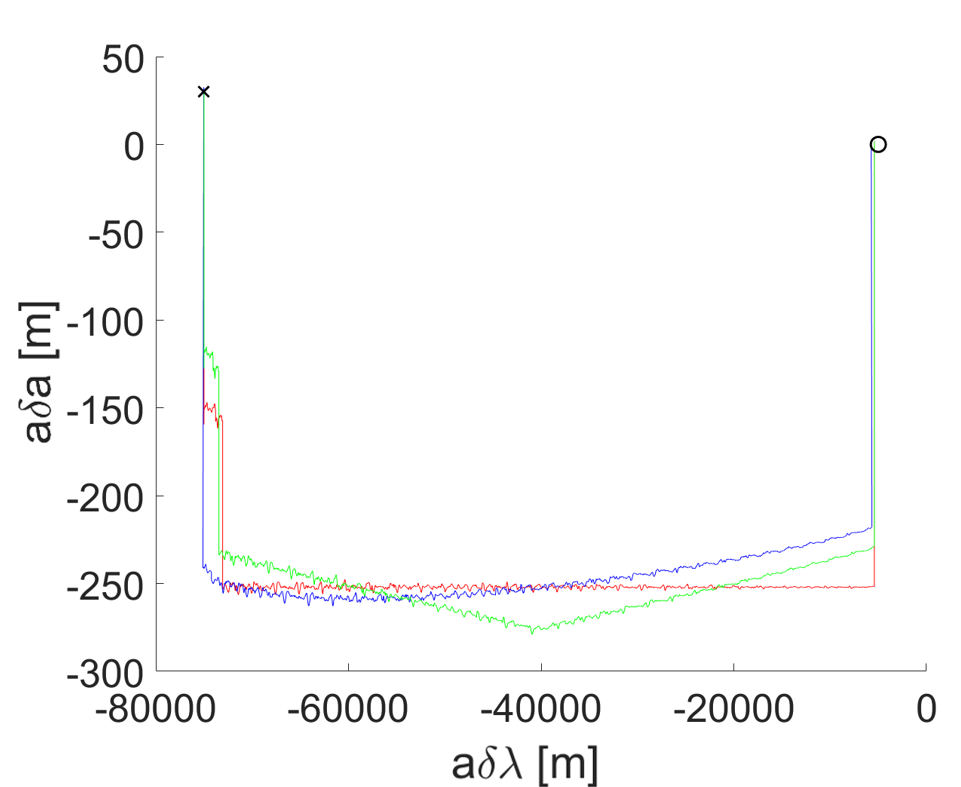}}
    \subfloat{\includegraphics[width=0.44\textwidth]{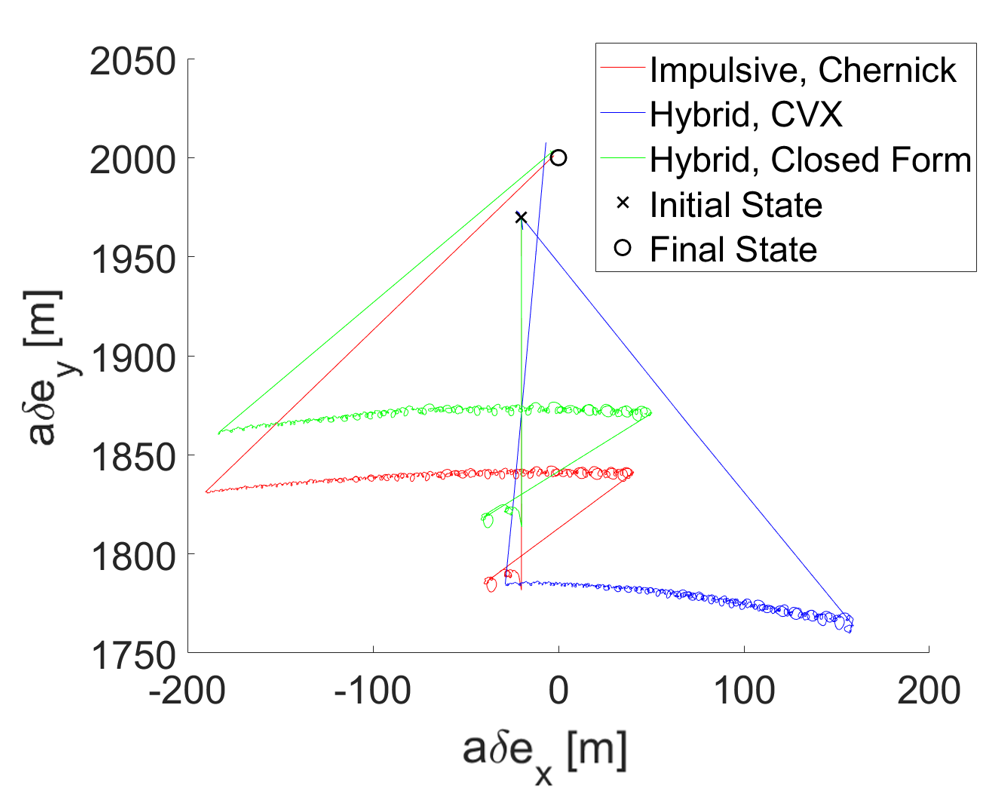}}
    \caption{In-plane ROE trajectories for Reconfiguration 2 under full-force dynamics model.}
    \label{fig: sim 2}
\end{figure}

Performance of the various control approaches is evaluated by comparing their final state errors. The type of reconfiguration that drives interest in hybrid control is a large change in along-track separation, corresponding to a $\delta \lambda$ dominant reconfiguration. Therefore, Reconfiguration 2 is chosen for this analysis. The in-plane ROE trajectories of the deputy spacecraft under Chernick impulsive control, numerically-solved hybrid control, and the proposed closed-form solutions for hybrid control are displayed in Figure \ref{fig: sim 2}, and the final state error achieved by each control method is listed in Table \ref{tab: sim 2 ROE error}.
\vspace{-5mm}
\begin{table}[h!]
    \centering
    \begin{tabular}{c|c|c|c|c|c|c}
    \hline \hline
    Control Approach & $a_c\delta a [m]$ & $a_c\delta\lambda [m]$ & $a_c\delta e_x [m]$ & $a_c\delta e_y [m]$ & $a_c\delta i_x [m]$ & $a_c\delta i_y [m]$ \\
    \hline
    Impulsive (Chernick) & 1.82 & -411.22 & -4.21 & 0.66 & 0.15 & -43.84 \\
    Hybrid (Numerical/CVX) & 2.78 & -735.04 & -6.88 & 7.51 & -3.09 & 0.85 \\
    Hybrid (Closed-form) & 1.82 & -402.96 & -4.49 & 3.39 & -1.37 & -48.33 \\
    \hline \hline
    \end{tabular}
    \caption{Final state error of control approaches for Reconfiguration 2 under full-force dynamics model.}
    \label{tab: sim 2 ROE error}
    \vspace{-4mm}
\end{table}

Given perfect knowledge of the spacecraft augmented ballistic coefficients, the proposed closed-form hybrid approach provides very similar control accuracy across all mean ROE to that of closed-form propulsive control. The trajectories are also very similar in Figure \ref{fig: sim 2}, with the effect of differential drag evident in the $\delta \boldsymbol{a}$ plane and the delta-v savings produced by hybrid control demonstrated in the $\delta \boldsymbol{e}$ plane by the size of the first maneuver from the initial state. The proposed hybrid approach demonstrates an unprecedented ability to accurately solve fully optimal roto-translational maneuver planning in closed-form. Furthermore, the numerical hybrid solution performs worse for in-plane ROE than either of the closed-form approaches, especially in the along-track or $\delta \lambda$ direction. The numerical solution also follows a very different trajectory in the $\delta \boldsymbol{e}$ plane. These differences occur because the numerical solver finds propulsive maneuvers with components in both the radial and tangential directions instead of just the tangential direction. These maneuvers have a larger effect on the $\delta \lambda$ pseudostate than tangential maneuvers alone, such that extra oscillations in the numerical solution for Reconfiguration 2 in Figure \ref{fig: drag profiles} are required to accommodate the $\delta \boldsymbol{e}$ pseudostate as well as the $\delta \lambda$ pseudostate to reach the same optimal in-plane cost.

The out-of-plane errors present in both the Chernick and proposed closed-form hybrid solutions are due to unmodelled changes in $\delta i_y$ that become significant over long control windows and are not a deficiency of the closed-form hybrid approach. These unmodelled changes primarily affect out-of-plane optimal propulsive maneuver timing, so inaccuracies could be mitigated by adding a heuristic to the control algorithm that favors optimal out-of-plane maneuvers occurring near the start of the control window. From the standpoint of passive safety, the reconfiguration still results in mostly parallel $\delta \boldsymbol{e}$ and $\delta \boldsymbol{i}$ vectors, so the loss in accuracy does not necessarily cause a concern for possible spacecraft collision. \cite{damico_phd}

\section{Conclusions and Future Work}

The future of Distributed Space Systems (DSS) will depend on increasingly optimal and accurate maneuvering to enhance the capabilities and mission lifetime of each small spacecraft. This work contributes directly to the state of the art by integrating propulsive and differential drag control to provide cost savings for large, time-constrained reconfigurations over propulsive-only methods. By posing the control problem in Relative Orbital Element (ROE) space, the relative dynamics of the spacecraft simplify to accurate linear models. Given that the effect of differential drag can be controlled on-board through spacecraft differential attitude, drag is incorporated into the control problem as a control force instead of a modification to the propagation of translational ROE state. The effect of a given differential drag profile can be subtracted from the original control problem to find a reduced remaining propulsive cost.

The application of reachable set theory produces optimality conditions for propulsive and drag control that are evaluated simultaneously to find hybrid delta-v minimum reconfiguration costs and optimal drag maneuver lengths. The reachable set itself is a geometric representation of the space of states that can be reached by a combination of maneuvers in a given control window at a specified total cost. The dominance case of a reconfiguration corresponds to the specific desired ROE that requires the highest cost to achieve, which drives the minimum cost of the overall reconfiguration. The simplicity of the reachable sets and corresponding limited number of reconfiguration dominance cases allows all globally optimal drag profiles to be derived analytically. A flexible drag profile framework allows portions of each of these profiles to be combined in closed-form under reasonable simplifying assumptions to create drag profiles for hybrid control completely decoupled from propulsive control computations. This approach also bridges the gap between hybrid and drag-only maneuver planning.

In general, hybrid control was shown to provide cost savings over propulsive-only methods. Specifically, the proposed algorithm provides roto-translational control solutions for hybrid control in an algebraic framework with high computational efficiency. These hybrid solutions were shown to provide provably global optimality for many reconfiguration types when compared with numerical approaches. The hybrid solutions also demonstrated the same final state accuracy performance as closed-form propulsive-only methods for full-force orbital simulations with relevant perturbations in Low Earth Orbit (LEO), given perfect state knowledge. The primary deficiency of hybrid control not covered by this paper is the drag control magnitude error due to uncertainty in the augmented differential ballistic coefficient over long control windows. This can be remedied with periodic control resolves throughout the control window, potentially with estimation updates to improve the characterization of the augmented ballistic coefficient of each spacecraft. Another potential flaw with the proposed approach is the assumption that propulsive and drag maneuvering can be conducted simultaneously, given that a spacecraft typically needs to change its attitude to point the on-board thrusters in any given direction. However, attitude changes to conduct propulsive maneuvers are relatively short in duration compared with the control windows relevant for hybrid control. Combined with the low thrust magnitude of differential drag, attitude maneuvering required for the propulsive maneuver plan can simply override the drag maneuver plan with limited impact on the final state error of the reconfiguration. The proposed algorithm's ability to solve optimal roto-translational control in closed-form will allow these advances to be implemented on a wide variety of spacecraft, including computationally-limited CubeSats and nanosats used on missions like SWARM-EX. The hybrid control approach will extend mission lifetimes and enhance reconfiguration capabilities for DSS in LEO.

\section{Acknowledgments}

This work was supported by the National Science Foundation (Award Number 1936512).

\bibliographystyle{AAS_publication}   
\bibliography{references.bib}   

\end{document}